\documentclass[aop,MSNbibl,seceqn,secthm,dvips]{arximspdf}
\makeatletter
   \@ifpackageloaded{graphicx}{}{\usepackage{graphicx}}
\makeatother
\usepackage{mathbh}

\doi{10.1214/15-AOP1017}
\volume{44}
\issue{3}
\pubyear{2016}
\firstpage{2024}
\lastpage{2063}
\docsubty{FLA}

\makeatletter
\newcommand{\E}{\mathbb{E}}
\newcommand{\FF}{\mathbb{F}}
\let\epsilon\varepsilon
\newtheorem{cor}[thm]{Corollary}
\newtheorem{lem}[thm]{Lemma}

\newproclaim{ex}[thm]{Example}
\newproclaim{ass}[thm]{Assumption}
\makeatother

\begin{document}
\begin{frontmatter}

\title{Berry--Esseen theorems under weak dependence}
\runtitle{Berry--Esseen theorems under weak dependence}

\begin{aug}
\author{\fnms{Moritz}~\snm{Jirak}\thanksref{T1}\corref{}\ead
[label=e1]{jirak@math.hu-berlin.de}}
\runauthor{M. Jirak}
\affiliation{Humboldt Universit\"{a}t zu Berlin}
\address{Institut f\"{u}r Mathematik\\
Unter den Linden 6\\
D-10099 Berlin\\
Germany\\
\printead{e1}}
\end{aug}
\thankstext{T1}{Supported by the Deutsche
Forschungsgemeinschaft via \textit{FOR 1735 Structural Inference in
Statistics: Adaptation and
Efficiency} is gratefully acknowledged.}

%
\received{\smonth{7} \syear{2014}}
%
\revised{\smonth{3} \syear{2015}}

%
\begin{abstract}
Let $ \{{X}_k \}_{k\geq\mathbb{Z}}$ be a stationary
sequence. Given
$p \in(2,3]$ moments and a mild weak dependence condition, we show a
Berry--Esseen theorem with optimal rate $n^{p/2 - 1}$. For $p \geq4$,
we also show a convergence rate of $n^{1/2}$ in $\mathcal{L}^q$-norm,
where $q \geq1$. Up to $\log n$ factors, we also obtain nonuniform
rates for any $p > 2$. This leads to new optimal results for many
linear and nonlinear processes from the time series literature, but
also includes examples from dynamical system theory. The proofs are
based on a hybrid method of characteristic functions, coupling and
conditioning arguments and ideal metrics.
\end{abstract}

%
\end{frontmatter}

\section{Introduction}\label{secintro}

Let $ \{X_k \}_{k \in\mathbb{Z}}$ be a zero mean process having
second moments $\mathbb{E} [X_k^2 ] < \infty$. Consider the partial
sum $S_n = \sum_{k = 1}^n X_k$ and its normalized variance $s_n^2 =
n^{-1}\operatorname{Var} [S_n ]$. A very important issue in
probability
theory and statistics is whether or not the central limit theorem
holds, that is, if we have
%
\begin{eqnarray}
\label{eqclt}
&& \lim_{n \to\infty}
 \Bigl|P \Bigl(S_n \leq x \sqrt
{n s_n^2} \Bigr) - \Phi(x) \Bigr| = 0,
\end{eqnarray}
where $\Phi(x)$ denotes the standard normal distribution function.
Going one step further, we can ask ourselves about the possible rate of
convergence in (\ref{eqclt}), more precisely, if it holds that
%
\begin{eqnarray}
\label{eqcltrate}
&& \lim_{n \to\infty}d (P_{S_n/\sqrt{n s_n^2}},
P_{Z} ) \mathfrak{r}_n < \infty\qquad \mbox{for a sequence }\mathfrak{r}_n \to \infty,
\end{eqnarray}
where $d (\cdot,\cdot )$ is a probability metric, $Z$
follows a standard normal distribution and~$P_{X}$ denotes the
probability measure induced by the random variable $X$. The rate
$\mathfrak{r}
_n$ can be considered as a measure of reliability for statistical
inference based on $S_n$, and large rates are naturally preferred. The
question of rate of convergence has been addressed under numerous
different setups with respect to the metric and underlying structure of
the sequence $ \{X_k \}_{k \in\mathbb{Z}}$ in the literature.
Perhaps one of the most important metrics is the Kolmogorov (uniform)
metric, given as
%
\begin{eqnarray}
\label{eqcltrateks}
&& \Delta_n = \sup_{x \in\mathbb{R}} \Bigl|P
\Bigl(S_n \leq x \sqrt {n s_n^2}
\Bigr) - \Phi(x) \Bigr|.
\end{eqnarray}
The latter has been studied extensively in the literature under many
different notions of (weak) dependence for $ \{X_k \}_{k \in
\mathbb{Z}}$. One general way to measure dependence is in terms of various
mixing conditions. In the case of the uniform metric, Bolthausen \cite{bolthausen1982ptrf} and Rio \cite{Rio1996} showed that it is
possible to obtain the rate $\mathfrak{r}_n = \sqrt{n}$ in (\ref
{eqcltrateks}), given certain mixing assumptions and a bounded
support of the underlying sequence $ \{X_k \}_{k \in\mathbb{Z}}$;
see also  \cite{borgnepene2005,chen2004aop,dedeckerprieur2005,hervepene2010bulletin}, among
others, for related results and extensions. Under the notion of $\alpha$-mixing, Tikhomirov \cite{tikhomirov1981} obtained $\mathfrak{r}_n =
n^{1/2}/(\log n)^2$, provided that $\mathbb{E} [|X_k|^3 ] <
\infty$
and the mixing coefficient decays exponentially fast; see also \cite{bentkus1997}. Martingales constitute another important class for the
study of (\ref{eqcltrateks}). Some relevant contributions in this
context are, for instance, Brown and Heyde \cite{brownheyde1970},
Bolthausen \cite{bolthausen1982} and more recently Dedecker et al.
\cite{dedeckerminimaldistance2009}. In the special case of
functionals of Gaussian or Poissonian sequences, deep results have been
obtained by Noudin and Peccati et al.; see, for instance, \cite{nourdin2009a,nourdin2009b} and \cite{peccati2010}. Another
stream of significant works focuses on stationary (causal)
Bernoulli-shift processes, given as
%
\begin{eqnarray}
\label{eqstructurecondition}
&& {X}_{k} = g_k (\epsilon_{k},
\epsilon_{k-1}, \ldots )\qquad \mbox{where }\{\epsilon_k
\}_{k \in\mathbb{Z}}\mbox{ is an i.i.d. sequence}.
\end{eqnarray}
The study of (\ref{eqcltrateks}) given the structure in (\ref
{eqstructurecondition}) has a long history, and dates back to Kac
\cite{kac1946} and Postnikov \cite{postnikov1966ergodic}.
Ibragimov \cite{ibramigov1967} established a rate of convergence,
$\mathfrak{r}_n = n^{1/2}/\sqrt{\log n}$, subject to an exponentially fast
decaying weak dependence coefficient. Using the technique of Tikhomirov
\cite{tikhomirov1981}, G\"{o}tze and Hipp obtained Edgeworth
expansions for processes of type (\ref{eqstructurecondition}) in a
series of works; cf. \cite{goetzehipp1983,goetzeHipp1989filedptrf,goetzehipp1995}; see also
Heinrich \cite{heinrich1990} and Lahiri \cite{lahiri1993}. This
approach, however, requires the validity of a number of technical
conditions. This includes in particular a conditional Cr\'{a}mer-like
condition subject to an exponential decay, which is somewhat difficult
to verify. In contrast, it turns out that a Berry--Esseen theorem only
requires a simple, yet fairly general dependence condition where no
exponential decay is required. Indeed, we will see that many popular
examples from the literature are within our framework. Unlike previous
results in the literature, we also obtain optimal rates for $p \in
(2,3)$ given (infinite) weak dependence, which to the best of our
knowledge is new (excluding special cases as linear processes). The
proofs are based on an $m$-dependent approximation ($m \to\infty$),
which is quite common in the literature. The substantial difference
here is the subsequent treatment of the $m$-dependent sequence. To
motivate one of the main ideas of the proofs, let us assume $p = 3$ for
a moment. Given a weakly, $m$-dependent sequence $ \{{X}_k \}
_{k \in\mathbb{Z}}$, one may show via classic arguments that
%
\begin{eqnarray}
\label{eqboundmom2.9ex}
&& \Delta_n \leq C\sqrt{m/n} \mathbb{E}
\bigl[|X_1|^3 \bigr],
\end{eqnarray}
provided that $\mathbb{E} [|X_1|^3 ] < \infty$ and $s_n^2 > 0$.
Note, however, since $X_k$ is weakly dependent, one finds that
%
\begin{eqnarray}
\label{eqboundmom3ex}
&& m^{-3/2} \bigl|\mathbb{E} \bigl[S_m^3
\bigr] \bigr| \leq\frac
{C}{\sqrt{m}}.
\end{eqnarray}
Hence if one succeeds in replacing $\mathbb{E} [|X_1|^3 ]$
in (\ref{eqboundmom2.9ex}) with (\ref{eqboundmom3ex}), one obtains the
optimal rate $\mathfrak{r}_n = \sqrt{n}$. A similar reasoning applies
to $p
\in(2,3)$. Unfortunately though, setting this idea to work leads to
rather intricate problems, and a technique like that of Tikhomirov
\cite{tikhomirov1981} is not fruitful, inevitably leading to a
suboptimal rate. Our approach is based on coupling and conditioning
arguments and ideal (Zolotarev) metrics. Interestingly, there is a
connection to more recent results of Dedecker et al. \cite{dedeckerminimaldistance2009}, who consider different (smoother)
probability metrics. We will see that at least some of the problems we
encounter may be redirected to these results after some preparation.

\section{Main results}\label{secmain}

Throughout\vspace*{1pt} this paper, we will use the following notation: for a random
variable $X$ and $p \geq1$, we denote with $ \|X \|_p =
\mathbb{E}
 [X^p ]^{1/p}$ the $\mathcal{L}^p$ norm. Let $ \{
\epsilon
_k \}_{k \in\mathbb{Z}}$ be a sequence of independent and identically
distributed random variables with values in a measurable space $\mathbb
{S}$. Denote the corresponding $\sigma$-algebra with $\mathcal{E}_k =
\sigma ( \epsilon_j, j \leq k  )$. Given a real-valued
stationary sequence $ \{{X}_k \}_{k\in\mathbb{Z}}$, we always
assume that $X_k$ is adapted to $\mathcal{E}_{k}$ for each $k \in
\mathbb{Z}$. Hence we implicitly assume that $X_{k}$ can be written as in (\ref
{eqstructurecondition}). For convenience, we write $X_k = g_k(\theta
_{k})$ with $\theta_k = (\epsilon_k, \epsilon_{k - 1},\ldots)$. The
class of processes that fits into this framework is large and contains
a variety of functionals of linear and nonlinear processes including
ARMA, GARCH and related processes (see, e.g., \cite{Gao2007,sipwu,tsay2005}), but also examples from dynamic system
theory. Some
popular examples are given below in Section~\ref{secex}. A nice
feature of the representation given in (\ref{eqstructurecondition})
is that it allows us to give simple, yet very efficient and general
dependence conditions. Following Wu \cite{wu2005}, let $ \{
\epsilon_k' \}_{k \in\mathbb{Z}}$ be an independent copy of
$ \{
\epsilon_k \}_{k \in\mathbb{Z}}$ on the same probability space, and
define the ``filter'' $\theta_{k}^{(l, \prime)}$ as
%
\begin{eqnarray}
\label{defnstrichdepe}
&& \theta_{k}^{(l,\prime)} = \bigl(\epsilon_{k},
\epsilon_{k - 1},\ldots,\epsilon_{k - l}',\epsilon_{k - l - 1},\ldots\bigr).
\end{eqnarray}
We put $\theta_{k}' = \theta_{k}^{(k, \prime)} = (\epsilon_{k}, \epsilon
_{k - 1},\ldots,\epsilon_0',\epsilon_{-1},\ldots)$ and ${X}_{k}^{(l,\prime)} =
g_k (\theta_{k}^{(l,\prime)}  )$, and in particular we set
${X}_{k}' ={X}_{k}^{(k,\prime)}$. As a dependence measure, we then consider
the quantity $\sup_{k \in\mathbb{Z}} \|{X}_{k} -
{X}_{k}^{(l,\prime)}
\|_p$, $p \geq1$. Dependence conditions of this type are quite general
and easy to verify in many cases; cf. \cite{aueetal2009,sipwu} and
the examples below. Observe that if the function $g = g_k$ does not
depend on $k$, we obtain the simpler version
%
\begin{eqnarray}
\label{eqhomodepcondi}
&& \sup_{k \in\mathbb{Z}} \bigl\|{X}_{k} -
{X}_{k}^{(l,\prime)} \bigr\|_p = \bigl\| {X}_{l} -
{X}_{l}^{\prime} \bigr\|_p.
\end{eqnarray}
Note that it is actually not trivial to construct a stationary process
$ \{X_k \}_{k \in\mathbb{Z}}$ that can only be represented
as $X_k
= g_k(\theta_k)$; that is, a function $g$ independent of $k$ such that
$X_k = g(\theta_k)$ for all $k \in\mathbb{Z}$ does not exist. We
refer to
Corollary~2.3 in Feldman and Rudolph \cite{feldmanrudolph1998} for
such an example.

We will derive all of our results under the following assumptions.

\begin{ass}\label{assmaindependence}
Let $ \{X_k \}_{k\in\mathbb{Z}}$ be stationary such that
for some
$p \geq2$:
\begin{longlist}[(iii)]
\item[(i)] $ \|X_k \|_p < \infty$, $\mathbb{E}
[X_k ] = 0$,
\item[(ii)] $\sum_{l = 1}^{\infty} l^2 \sup_{k \in\mathbb
{Z}} \|
{X}_{k} - {X}_{k}^{(l,\prime)}  \|_p < \infty$,
\item[(iii)] $s^2 > 0$, where $s^2 = \sum_{k \in\mathbb{Z}}\mathbb
{E}
[X_0X_k ]$.
\end{longlist}
\end{ass}

In the sequel, ${B}$ denotes a varying absolute constant, depending
only on $p$, $\sum_{l = 1}^{\infty} l^2 \sup_{k \in\mathbb
{Z}} \|
{X}_{k} - {X}_{k}^{(l,\prime)}  \|_p$ and $s^2$. The following
theorem is one of the main results of this paper.
%
\begin{thm}\label{thmberry}
Grant Assumption~\ref{assmaindependence} for some $p \in(2,3]$, and
let $s_n^2 = n^{-1} \|S_n  \|_2^2$. Then
\begin{eqnarray*}
&& \sup_{x \in\mathbb{R}} \Bigl|P \Bigl(S_n\big/\sqrt{n
s_n^2} \leq x \Bigr) - \Phi (x ) \Bigr| \leq\frac{{B}}{n^{p/2 - 1}},
\end{eqnarray*}
and hence we may select $\mathfrak{r}_n = n^{p/2 - 1}$.
\end{thm}

Theorem~\ref{thmberry} provides optimal convergence rates under mild
conditions. In particular, it seems that this is the first time optimal
rates are shown to hold under general infinite weak dependence
conditions if $p \in(2,3)$. Examples to demonstrate the versatility of
the result are given in Section~\ref{secex}. In particular, we
consider functions of the dynamical system $Tx = 2x \operatorname{mod}1$ in Example~\ref{exnumber}, a problem which has been studied in
the literature
for decades. Combining Theorem~\ref{thmberry} with results of
Dedecker and Rio \cite{dedeckerriomean2008}, we also obtain
optimal results for the ${\mathcal L}^q$-norm for martingale differences.
%
\begin{thm}\label{thmLp}
Grant Assumption~\ref{assmaindependence} for some $p \geq4$, and
let $s_n^2 = n^{-1} \|S_n  \|_2^2$. If $ \{X_k \}_{k\in\mathbb{Z}}$ is a martingale difference sequence, then for any $q \geq1$ we have
\begin{eqnarray*}
&& \int_{\mathbb{R}} \Bigl|P \Bigl(S_n\big/\sqrt{n
s_n^2} \leq x \Bigr) - \Phi (x ) \Bigr|^q \,dx \leq
{B}n^{-q/2}.
\end{eqnarray*}
\end{thm}
Note that in the case $q = 1$, the results of Dedecker and Rio
\cite{dedeckerriomean2008} are more general. The nonuniform analogue to
Theorems \ref{thmberry} and \ref{thmLp} is given below. Here, we
obtain optimality up to logarithmic factors.
%
\begin{thm}\label{thmnonunif}
Grant Assumption~\ref{assmaindependence} for some $p > 2$. Then for
any $x \in\mathbb{R}$,
\begin{eqnarray*}
&& \Bigl|P \Bigl(S_n\big/\sqrt{n s_n^2} \leq x
\Bigr) - \Phi (x ) \Bigr| \leq n^{-(p \wedge 3)/2 + 1} \frac{{B} ( \log n
)^{p/2}}{1 + |x|^{p}},
\end{eqnarray*}
where $a\wedge b = \min\{a,b\}$.
\end{thm}

As a particular application of Theorem~\ref{thmnonunif}, consider
$f (|S_n|/\sqrt{n s_n^2} )$ where the function $f(\cdot)$ satisfies
%
\begin{eqnarray}
\label{eqgconditforapprox}
&& f (0 ) = 0 \quad \mbox{and}\quad \int_0^{\infty}
\frac
{|f'(x)|}{1 + |x|^{p}} \,dx < \infty
\end{eqnarray}
for some $p > 0$, and the derivative $f'(x)$ exists for $x \in
(0,\infty)$. If $ \|S_n \|_p < \infty$, property (\ref
{eqgconditforapprox}) implies the identity
\begin{eqnarray*}
&& \mathbb{E} \Bigl[f\Bigl(|S_n|\big/\sqrt{n s_n^2}
\Bigr) \Bigr] = \int_{0}^{\infty} f'(x)
P \Bigl(|S_n|\big/\sqrt{n s_n^2} \geq x
\Bigr) \,dx,
\end{eqnarray*}
and we thus obtain the following corollary.

\begin{cor}\label{corgconvergence}
Grant Assumption~\ref{assmaindependence} for some $p > 2$. If (\ref
{eqgconditforapprox}) holds, then
\begin{eqnarray*}
&& \biggl|\mathbb{E} \Bigl[f\Bigl(|S_n|\big/\sqrt{n s_n^2}
\Bigr) \Bigr] - \int_{\mathbb{R}}f \bigl(|x| \bigr)\,d \Phi(x) \biggr| \leq
{B}n^{-(p \wedge 3)/2 + 1} ( \log n )^{p/2}.
\end{eqnarray*}
\end{cor}

As a special case, consider $f(|x|) = |x|^q$, $q > 0$. We may then use
Corollary~\ref{corgconvergence} to obtain rates of convergence for moments.

\begin{cor}\label{cormomconvergence}
Grant Assumption~\ref{assmaindependence} for some $p > 2$. Then for
any $0 < q < p$, we have
\begin{eqnarray*}
&& \biggl| \Bigl\|S_n\big/\sqrt{n s_n^2}
\Bigr\|_q^q - \int_{\mathbb{R}} |x|^q
\,d\Phi(x) \biggr| \leq{B}n^{-(p \wedge 3)/2 + 1} (\log n )^{p/2}.
\end{eqnarray*}
\end{cor}

In the special case of i.i.d. sequences and $0 < p < 4$, sharp results
in this context have been obtained in Hall \cite{hall1982}. It seems
that related results for dependent sequences are unknown.

\section{Applications and examples}\label{secex}

All examples considered here are time-homogenous Bernoulli-shift
processes; that is, $g = g_k$ does not depend on $k$, and hence
equality (\ref{eqhomodepcondi}) holds.

\begin{ex}[(Functions of linear process)]
Let $\mathbb{S} = \mathbb{R}$, and suppose that the sequence $ \{
\alpha
_i \}_{i \in\mathbb{N}}$ satisfies $\sum_{i = 0}^{\infty}
\alpha_i^2
< \infty$. If $\|\epsilon_k\|_2 < \infty$, then one may show that
the linear process
\begin{eqnarray*}
&& Y_k = \sum_{i = 0}^{\infty}
\alpha_i \epsilon_{k - i}\qquad \mbox{exists and is stationary}.
\end{eqnarray*}
Let $f$ be a measurable function such that $\mathbb{E} [X_k
] = 0$,
where $X_k = f(Y_k)$. If $f$ is H\"{o}lder continuous with regularity
$0 < \beta\leq1$, that is, $ |f(x) - f(y) | \leq c
|x-y|^{\beta}$, then for any $p \geq1$
\begin{eqnarray*}
&& \bigl\|X_k - X_k' \bigr\|_p \leq c
\alpha_k^{\beta} \|\epsilon _0 \|_p.
\end{eqnarray*}
Hence if $\sum_{i = 0}^{\infty} i^2 |\alpha_i|^{\beta} < \infty$
and $s^2 > 0$, then Assumption~\ref{assmaindependence} holds.
\end{ex}

\begin{ex}[{[Sums of the form $\sum f (t 2^k )$]}]\label{exnumber}
Consider the measure preserving transformation $Tx = 2x \operatorname{mod}1$ on the
probability space $ ([0,1],\mathcal{B},\bolds{\lambda} )$, with
Borel $\sigma$-algebra $\mathcal{B}$ and Lebesgue measure $\bolds{\lambda}
$. Let $U_0 \thicksim  \operatorname{Uniform}[0,1]$. Then $T U_0 = \sum_{j =
0}^{\infty} 2^{-j-1}\zeta_j$, where $\zeta_j$ are Bernoulli random
variables. The flow $T^k U_0$ can then be written as $T^k U_0 = \sum_{j = 0}^{\infty} 2^{-j-1} \zeta_{j + k}$; see
\cite{ibramigov1967}. The study about the behavior of $S_n = \sum_{k =
1}^{n} f (T^k U_0 )$ for appropriate functions $f$ has a very
long history and dates back to Kac \cite{kac1946}. Since then,
numerous contributions have been made; see, for instance,
\cite{berkes2014,billingsley1999,denkerkeller1986,dedeckerriomean2008,hormann2009,ibramigov1967,ladokhin1971,mcleish1975,moskvin1979local,petit1992,postnikov1966ergodic},
to name a few. Here, we consider the
following class of functions. Let $f$ be a function defined on the unit
interval $[0,1]$, such that
%
\begin{eqnarray}
&& \int_0^1 f(t) \,d t =
0, \qquad \int_0^1 \bigl|f(t)\bigr|^p \,d t < \infty
\quad\mbox{and}
\nonumber
\\[-8pt]
\label{eqibracondi}
\\[-8pt]
\nonumber
&&\int_0^1 t^{-1} \bigl|
\log(t)\bigr|^2 w_p (f,t )\, dt < \infty,
\end{eqnarray}
where $w_p (f,t )$ denotes a $\mathcal{L}^p
([0,1], \bolds{\lambda} )$ modulos of continuity of $f \in
\mathcal
{L}^p ([0,1],\break\bolds{\lambda} )$. This setup is a little more
general than in \cite{ibramigov1967}. For $x \in\mathbb{R}^+$, let
$\bar
{f}(x) = f(x - \lfloor x \rfloor)$; that is, $\bar{f}$ is the
one-periodic extension\vspace*{1pt} to the positive real line. One then often finds
the equivalent formulation $S_n = \sum_{k = 1}^n \bar{f} (2^k
U_0 )$ in the literature.\vspace*{1pt} Consider now the partial sum $S_n = \sum_{k = 1}^n \bar{f} (2^k U_0 )$. Ibragimov
\cite{ibramigov1967} showed that
%
\begin{eqnarray}
\label{eqexnumber1}
&& \sup_{x \in\mathbb{R}} \Bigl|P \Bigl(S_n\big/\sqrt
{n s_n^2} \leq x \Bigr) - \Phi (x ) \Bigr| \leq C
\biggl(\frac{\log n}{n} \biggr)^{p/2-1}.
\end{eqnarray}
%
By alternative methods, according to ~\cite{dedeckerriomean2008}, the results of ~\cite{jan:2001} allow to remove
the logarithmic factor if $\|f\|_{\infty} < \infty$. A priori, the
sequence $ \{T^k U_0 \}_{k \in\mathbb{Z}}$ does not
directly fit
into our framework, which, however, can be achieved by a simple time
flip. Define the function $T_n(i) = n - i + 1$ for $i \in\{n,n-1,\ldots\}
$, and let $\epsilon_k = \zeta_{T_n(k)}$. Then we may write
\begin{eqnarray*}
&& X_k = f \bigl(T^k U_0 \bigr) = f \Biggl(
\sum_{j = 0}^{\infty} \epsilon_{k - j}
2^{-j-1} \Biggr), \qquad k \in \{1,\ldots ,n \}.
\end{eqnarray*}
Note that we have to perform this time flip for every $n \in\mathbb{N}$,
which, however, has no impact on the applicability of our results.
Using the same arguments as in \cite{ibramigov1967}, we find that
(\ref{eqibracondi}) implies that for $p \in(2,3]$
\begin{eqnarray*}
&& \sum_{k = 1}^{\infty} k^2
\bigl\|X_k - X_k' \bigr\|_p < \infty.
\end{eqnarray*}
If $s^2 > 0$, we see that Assumption~\ref{assmaindependence}
holds. In particular, an application of Theorem~\ref{thmberry} gives
the rate $\mathfrak{r}_n = n^{p/2 - 1}$, thereby removing the
unnecessary $\log
n$ factor in (\ref{eqexnumber1}) for the whole range $p \in(2,3]$.
\end{ex}

\begin{ex}[($m$-dependent processes)]\label{exmdep}
$\!\!$Consider the zero mean $m$-\break dependent process $Y_k = f (\zeta
_k,\ldots,\zeta_{k-m+1} )$, where $m \in\mathbb{N}$ and $f$ is a
measurable function and $ \{\zeta_k \}_{k \in\mathbb{Z}}$
is i.i.d.
and takes values in $\mathbb{S}$. $m$ may depend on $n$ such that $n/m
\to\infty$, but we demand in addition that
%
\begin{eqnarray}
\label{eqmdepexvarcondi}
&& \mathop{\liminf}_{n \to\infty}\operatorname{Var} \Biggl[\sum
_{k = 1}^n Y_k \Biggr]\Big/ (n m ) >
0.
\end{eqnarray}
In this context, it is useful to work with the transformed block-variables
\begin{eqnarray*}
&& X_k = \frac{1}{m}\sum_{l = 0}^{m-1}Y_{mk-l},
\qquad k \in\mathbb{Z},
\end{eqnarray*}
and write $X_k = g (\epsilon_k,\epsilon_{k-1} )$ where
$\epsilon_k =  (\zeta_{km},\ldots,\zeta_{(k-1)m + 1}
)^{\top} \in\mathbb{S}^m$; hence $ \{X_k \}_{k \in\mathbb{Z}}$
is a two-dependent sequence. This representation ensures that
Assumption~\ref{assmaindependence}(i) and (ii) hold for $ \{
X_k \}_{k \in\mathbb{Z}}$, independently of the value of $m$. The
drawback of this block-structure is that we loose a factor $m$, since
we have
\begin{eqnarray*}
&& \frac{1}{\sqrt{n m}}S_n = \frac{1}{\sqrt{n m}}\sum
_{k = 1}^n Y_k = \frac{1}{\sqrt{n/m}}\sum
_{k = 1}^{n/m} X_k,
\end{eqnarray*}
where we assume that $n/m \in\mathbb{N}$ for simplicity. However,
this loss
is known in the literature: Theorem~\ref{thmberry} now yields the
commonly observed rate $\mathfrak{r}_n = (n/m)^{p/2-1}$ in the context of
$m$-dependent sequences satisfying\vspace*{1pt} (\ref{eqmdepexvarcondi}); see,
for instance, Theorem~2.6 in \cite{chen2004aop}. In the latter, the
rate $\mathfrak{r}_n = (n/m)^{p/2-1}$ is not immediately obvious, but follows
from elementary computations using (\ref{eqmdepexvarcondi}).
\end{ex}

\begin{ex}[(Iterated random function)]
Iterated random functions (cf. \cite{diaconisfreedman1999}) are an
important class of processes. Many nonlinear models like ARCH, bilinear
and threshold autoregressive models fit into this framework. Let
$\mathbb{S} = \mathbb{R}$ and $ \{X_k \}_{k \in\mathbb
{Z}}$ be defined
via the recursion
\begin{eqnarray*}
&& X_k = G (X_{k-1},\epsilon_k ),
\end{eqnarray*}
commonly referred to as \textit{iterated random functions}; see, for
instance, \cite{diaconisfreedman1999}. Let
%
\begin{eqnarray}
&& L_{\epsilon} = \sup_{x \neq y} \frac{ |G(x,\epsilon) -
G(y,\epsilon) |}{|x-y|}
\end{eqnarray}
be the Lipschitz coefficient. If $ \|L_{\epsilon} \|_p < 1$
and $ \|G(x_0,\epsilon) \|_p < \infty$ for some $x_0$, then
$X_k$ can be represented as $X_k = g (\epsilon_k,\epsilon
_{k-1},\ldots )$ for some measurable function $g$. In addition, we have
%
\begin{eqnarray}
&& \bigl\|X_k - X_k' \bigr\|_p \leq C
\rho^{-k}\qquad  \mbox{where }0 < \rho< 1;
\end{eqnarray}
see \cite{wushaoiterated2004}. Hence if $\mathbb{E}[X_k] = 0$ and $s^2
> 0$, Assumption~\ref{assmaindependence} holds. As an example,
consider the \textit{stochastic recursion}
\begin{eqnarray*}
&& X_{k+1} = a_{k+1} X_{k} + b_{k+1},\qquad k
\in\mathbb{Z},
\end{eqnarray*}
where $ \{a_k,b_k \}_{k \in\mathbb{Z}}$ is an i.i.d.
sequence. Let
$\epsilon_k = (a_k,b_k)$. If we have, for some $p \geq2$,
%
\begin{eqnarray}
\label{stocheqcondi}
&& \|a_k\|_p < 1\quad \mbox{and}\quad
\|b_k\|_p < \infty,\qquad \mathbb{E} [b_k ] = 0,
\end{eqnarray}
then $ \|L_{\epsilon} \|_p \leq\|a_k\|_p < 1$, and
Assumption~\ref{assmaindependence} holds if $s^2 > 0$. In
particular, if $a_k,b_k$ are independent, then one readily verifies that
\begin{eqnarray*}
&& s^2 = \frac{\|b_0\|_2^2}{1 - \|a_0\|_2^2} \biggl(1 + \frac{2\mathbb{E}
[a_0]}{1 - \mathbb{E}[a_0]} \biggr),
\end{eqnarray*}
which is strictly positive since $ |\mathbb{E}[a_0] | < 1$ by
Jensen's inequality. Hence if (\ref{stocheqcondi}) holds for $p >
2$, then Assumption~\ref{assmaindependence} holds for $p$. Analogue
conditions can be derived for higher order recursions.
\end{ex}

\begin{ex}[{[$\operatorname{GARCH}(\mathfrak{p},\mathfrak{q})$ sequences]}]\label{exgarch}
Let $\mathbb{S} = \mathbb{R}$. Another very prominent stochastic
recursion is
the $\operatorname{GARCH}(\mathfrak{p},\mathfrak{q})$ sequence, given through the relations
\begin{eqnarray*}
X_k &=& \epsilon_k L_{k} \qquad \mbox{where }\{\epsilon_{k} \}_{k \in\mathbb{Z}}\mbox{ is a zero mean i.i.d. sequence and}
\\
L_k^2  &=&  \mu+ \alpha_1 L_{k - 1}^2
+ \cdots + \alpha_\mathfrak{p} L_{k
- \mathfrak{p}}^2 +
\beta_1 X_{k - 1}^2 + \cdots +
\beta_{\mathfrak
{q}} X_{k - \mathfrak{q}}^2,
\end{eqnarray*}
with $\mu, \alpha_1,\ldots,\alpha_\mathfrak{p}, \beta_1,\ldots,\beta
_{\mathfrak{q}} \in\mathbb{R}$. We assume that $ \|\epsilon
_k \|
_p < \infty$ for some $p \geq2$. An important quantity is
\begin{eqnarray*}
&& \gamma_C = \sum_{i = 1}^{r}
\bigl\|\alpha_i + \beta_i \epsilon _i^2
\bigr\|_2 \qquad \mbox{with $r = \max\{\mathfrak{p},\mathfrak {q}\}$},
\end{eqnarray*}
where we replace possible undefined $\alpha_i, \beta_i$ with zero. If
$\gamma_C < 1$, then $ \{X_k \}_{k \in\mathbb{Z}}$ is stationary;
cf. \cite{bougerol1992aop}. In particular, it was shown in \cite{berkes2008letter} that $ \{X_k \}_{k \in\mathbb{Z}}$ may be
represented as
\begin{eqnarray*}
&& X_k = \sqrt{\mu}\epsilon_k \Biggl(1 + \sum
_{n = 1}^{\infty} \sum_{1 \leq l_1,\ldots,l_n\leq r}
\prod_{i = 1}^n \bigl(\alpha_{l_i}
+ \beta_{l_i}\epsilon_{j - l_1 - \cdots - l_i}^2 \bigr)
\Biggr)^{1/2}.
\end{eqnarray*}
Using this representation and the fact that $|x-y|^p \leq|x^2 -
y^2|^{p/2}$ for $x,y \geq0$, $p \geq1$, one can follow the proof of
Theorem~4.2 in \cite{aueetal2009} to show that
\begin{eqnarray*}
&& \bigl\|X_k - X_k' \bigr\|_p \leq C
\rho^{k}\qquad  \mbox{where $0 < \rho< 1$.}
\end{eqnarray*}
Since $\mathbb{E} [X_k ] = \mathbb{E} [\epsilon
_k ] = 0$,
Assumption~\ref{assmaindependence} holds if $s^2 > 0$. We remark
that previous results on $\Delta_n$, in the case of GARCH($\mathfrak
{p},\mathfrak{q})$ sequences, either require heavy additional
assumptions or have suboptimal rates; cf. \cite{hormann2009}.
\end{ex}

\begin{ex}[(Volterra processes)]
In the study of nonlinear processes, Volterra processes are of
fundamental importance. Following Berkes et al. \cite{berkes2014}, we consider
\begin{eqnarray*}
&& X_k = \sum_{i = 1}^{\infty} \sum
_{0 \leq j_1 <\cdots<j_i} a_k (j_1,\ldots,
j_i ) \epsilon_{k - j_1} \cdots\epsilon_{k-j_i},
\end{eqnarray*}
where $\mathbb{S} = \mathbb{R}$ and $ \|\epsilon_k \|_p <
\infty$
for $p \geq2$, and $a_k$ are called the $k$th Volterra kernel. Let
\begin{eqnarray*}
&& A_{k,i} = \sum_{k \in\{j_1,\ldots,j_i\}, 0 \leq j_1 <\cdots<j_i} \bigl|a_k
(j_1,\ldots, j_i ) \bigr|.
\end{eqnarray*}
Then there exists a constant $C$ such that
\begin{eqnarray*}
&& \bigl\|X_k - X_k' \bigr\|_p \leq C\sum
_{i = 1}^{\infty} \| \epsilon_0
\|_p^{i} A_{k,i}.
\end{eqnarray*}
Thus if $\sum_{k,i = 1}^{\infty} k^2 A_{k,i} < \infty$ and $s^2 >
0$, then Assumption~\ref{assmaindependence} holds.
\end{ex}

\section{Proofs}\label{secproofs}

The main approach consists of an $m$-dependent approximation where $m
\to\infty$, followed by characteristic functions and Esseen's
inequality. However, here the trouble starts, since we cannot factor
the characteristic function as in the classic proof, due to the
$m$-dependence. Tikhomirov \cite{tikhomirov1981} uses a
chaining-type argument, which is also fruitful for Edgeworth
expansions; cf. \cite{goetzehipp1983}. However, since this approach
inevitably leads to a loss in the rate, this is not an option for
Berry--Esseen-type results. In order to circumvent this problem, we
first work under an appropriately chosen conditional probability
measure $P_{\mathbb{F}_m}$. Unfortunately though, this leads to rather
intricate problems, since all involved quantities of interest are then
random. We first consider the case of a weakly $m$-dependent sequence
$ \{X_k \}_{k \in\mathbb{Z}}$, where $m \to\infty$ as $n$
increases. Note that this is different from Example~\ref{exmdep}.
For the general case, we then construct a suitable $m$-dependent
approximating sequence such that the error of approximation is
negligible, which is carried out in Section~\ref{secproofofmaintheorem}. The overall proof of Theorem~\ref
{thmberry} is lengthy. Important technical auxiliary results are
therefore established separately in Section~\ref{secmainlemmas}.
Minor additionally required results are collected in Section~\ref{secauxresults}. The proofs of Theorems~\ref{thmLp} and \ref
{thmnonunif} are given in Sections~\ref{secprooflp} and \ref
{secthmproofnonunif}. To simplify the notation in the proofs, we
restrict ourselves to the case of homogeneous Bernoulli shifts, that
is, where $X_k = g (\epsilon_k,\epsilon_{k-1},\ldots )$,
and the function $g$ does not depend on $k$. This requires
substantially fewer indices and notation throughout the proofs, and, in
particular, (\ref{eqhomodepcondi}) holds. The more general
nonhomogenous (but still stationary) case follows from straightforward
(notational) adaptations. This is because the key ingredient we require
for the proof is the Bernoulli-shift structure (\ref
{eqstructurecondition}) in connection with the summability condition,
Assumption~\ref{assmaindependence}(ii). Whether or not $g$ depends
on $k$ is of no relevance in this context.
%

\subsection{$m$-dependencies}

In order to deal with $m$-dependent sequences, we require some
additional notation and definitions. Throughout the remainder of this
section, we let
\begin{eqnarray*}
&& X_k = f_m(\epsilon_k,\ldots,
\epsilon_{k-m+1})\qquad \mbox{for $m \in\mathbb{N}$, $k \in\mathbb{Z}$,}
\end{eqnarray*}
and measurable functions $f_m \dvtx  \mathbb{S}^{m} \to \mathbb{R}$,
where $m =
m_n \to\infty$ as $n$ increases. We work under the following conditions:

\begin{ass}\label{assdependence}
Let $ \{X_k \}_{k\geq\mathbb{Z}}$ be such that for some $p
\geq2$,
uniformly in $m$:
\begin{longlist}[(iii)]
\item[(i)] $ \|X_k \|_p < \infty$, $\mathbb{E}
[X_k ] = 0$,\vspace*{1pt}
\item[(ii)] $\sum_{k = 1}^{\infty} k^2  \|X_k - X_k' \|_p
< \infty$,\vspace*{1pt}
\item[(iii)] $s_m^2 > 0$,
\end{longlist}
where $s_m^2 = \sum_{k \in\mathbb{Z}}\mathbb{E} [X_0X_k ]
= \sum_{k =
-m}^{m}\mathbb{E} [X_0X_k ]$.
\end{ass}

Observe that this setup is fundamentally different from that considered
in Example~\ref{exmdep}. In particular, here we have that
$\operatorname{Var}
 [S_n ] \thicksim n $. Define the following $\sigma$-algebra:
%
\begin{eqnarray}
\label{defnsigmaalgebra}
&& \mathbb{F}_m = \sigma \bigl(\epsilon_{-m+1},
\ldots,\epsilon_0,\epsilon _1',\ldots,
\epsilon_m',\epsilon_{m+1},\ldots,
\epsilon_{2m},\epsilon _{2m+1}',\ldots \bigr),
\end{eqnarray}
where we recall that $ \{\epsilon_k \}_{k \in\mathbb{Z}}$ and
$ \{\epsilon_k' \}_{k \in\mathbb{Z}}$ are mutually independent,
identically distributed random sequences. We write $P_{\mathbb
{F}_m}(\cdot)$
for the conditional law and $\mathbb{E}_{\mathbb{F}_m}[\cdot]$ (or
$\mathbb{E}_{\mathcal
{H}}[\cdot]$) for the conditional expectation with respect to $\mathbb{F}_m$
(or some other $\sigma$-algebra~$\mathcal{H}$). We introduce
\begin{eqnarray*}
&& S_{|m}^{(1)} = \sum_{k = 1}^n
X_k - \mathbb{E} [X_k|\mathbb {F}_m ]
\quad \mbox{and}\quad  S_{|m}^{(2)} = \sum_{k = 1}^n
\mathbb{E} [X_k|\mathbb{F} _m ],
\end{eqnarray*}
hence
\begin{eqnarray*}
&& S_n = \sum_{k = 1}^n
X_k = S_{|m}^{(1)} + S_{|m}^{(2)}.
\end{eqnarray*}
To avoid any notational problems, we put $X_k = 0$ for $k \notin
 \{1,\ldots,n \}$. Let $n = 2(N-1)m + m'$, where $N,m$ are
chosen such that $c_0 m \leq m' \leq m$ and $c_0 > 0$ is an absolute
constant, independent of $m,n$. For $1 \leq j \leq N$, we construct the
block random variables
\begin{eqnarray*}
&& U_j = \sum_{k = (2j-2)m + 1}^{(2j - 1)m}
X_k - \mathbb{E} [X_k|\mathbb{F} _m ]
\quad \mbox{and}\quad R_j = \sum_{k = (2j-1)m + 1}^{2jm}
X_k - \mathbb{E} [X_k|\mathbb{F}_m ],
\end{eqnarray*}
and put $Y_j^{(1)} = U_j + R_j$, hence $S_{|m}^{(1)} = \sum_{j =
1}^{N} Y_j^{(1)}$. Note that by construction of the blocks,
$Y_j^{(1)}$, $j = 1,\ldots,N$ are independent random variables under the
conditional probability measure $P_{\mathbb{F}_m}(\cdot)$, and are
identically distributed at least for $j = 1,\ldots,N-1$ under $P$.
%
%
We also put $Y_1^{(2)} = \sum_{k = 1}^{m} \mathbb{E} [X_k
|\mathbb{F}
_m ]$ and $Y_j^{(2)} = \sum_{k = (j-1)m + 1}^{(j+1)m} \mathbb
{E}
[X_k |\mathbb{F}_m ]$ for $j = 2,\ldots,N$. Note\vspace*{1pt} that
$Y_j^{(2)}$, $j
= 1,\ldots,N$ is a sequence of independent random variables. The following
partial and conditional variances are relevant for the proofs:
\begin{eqnarray*}
\sigma_{j|m}^2 &=& \frac{1}{2m}
\mathbb{E}_{\mathbb{F}_m} \bigl[ \bigl(Y_j^{(1)}
\bigr)^2 \bigr]\quad \mbox{and}\quad \sigma_j^2 =
\mathbb{E} \bigl[\sigma_{j|m}^2 \bigr],
\\
\sigma_{|m}^2 &=& \frac{1}{n}\mathbb{E}
\bigl[\bigl(S_{|m}^{(1)}\bigr)^2 |
\mathbb{F}_m \bigr] = \frac{1}{N + m'/2m} \sum
_{j = 1}^N \sigma_{j|m}^2,
\\
\overline{\sigma}_m^2 &=& \mathbb{E} \bigl[
\sigma_{|m}^2 \bigr] = \frac
{1}{N + m'/2m} \sum
_{j = 1}^N \sigma_j^2,
\\
\widehat{\sigma}_m^2 &=& \frac{1}{2m}\sum
_{k = 1}^m\sum_{l = 1}^m
\mathbb{E} [X_kX_l ].
\end{eqnarray*}
As we shall see below, these quantities are all closely connected. Note
that $\sigma_i^2 = \sigma_j^2$ for $1 \leq i,j \leq N-1$, but $\sigma
_1^2 \neq\sigma_{N}^2$ in general. Moreover, we have the equation
\begin{eqnarray}
\label{decompsigma}
&& 2m \widehat{\sigma}_m^2 = m
s_m^2 - \sum_{k \in\mathbb{Z}} m \wedge
|k|\mathbb{E} [X_0 X_k ].
\end{eqnarray}
The above relation is important, since Lemma~\ref{lemsigexpressionsrelations} yields that under Assumption~\ref
{assdependence} we have $2\widehat{\sigma}_m^2 = s_m^2 + \mathcal{O}
 (m^{-1} )$. Moreover, Lemma~\ref{lemsigexpansion} gives
$\sigma_j^2 = \widehat{\sigma}_m^2 + \mathcal{O} (m^{-1}
)$ for
$1 \leq j \leq N-1$. We conclude that
%
\begin{eqnarray}
\label{eqsigmabiggerzero}
&& \sigma_j^2 = s_m^2/2
+ \mathcal{O} \bigl(m^{-1} \bigr) > 0\qquad \mbox{for sufficiently large
$m$.}
\end{eqnarray}
The same is true for $\sigma_N^2$, since $m' \geq c_0 m$. Summarizing,
we see that we do not have any degeneracy problems for the partial
variances $\sigma_j^2$, $1 \leq j \leq N$ under Assumption~\ref
{assdependence}. For the second part $S_{|m}^{(2)}$, we introduce
$\overline{\varsigma}_{m}^2 = n^{-1} \|S_{|m}^{(2)} \|
_2^2$. One then readily derives via conditioning arguments that
%
\begin{eqnarray}
\label{defnsnm}
&& s_{nm}^2 \stackrel{\mathrm{def}}{=}
n^{-1} \|S_n \|_2^2 =
n^{-1} \bigl\|S_{|m}^{(1)} \bigr\|_2^2
+ n^{-1} \bigl\|S_{|m}^{(2)} \bigr\| _2^2
= \overline{\sigma}_{m}^2 + \overline{
\varsigma}_{m}^2.
\end{eqnarray}

We are now ready to give the main result of this section.
%
\begin{thm}\label{thmmdependent}
Grant Assumption~\ref{assdependence}, and let $p \in(2,3]$. Assume
in addition that $N = N_n = n^{\lambda}$ for $0 < \lambda\leq p/(2p +
2)$. Then
\begin{eqnarray*}
&& \sup_{x \in\mathbb{R}} \bigl|P (S_n/\sqrt{n} \leq x ) - \Phi
(x/s_{nm} ) \bigr| \leq c(\lambda,p) n^{-p/2 + 1},
\end{eqnarray*}
where $c(\lambda,p)>0$ depends on $\lambda$, $p$, $\sum_{k =
1}^{\infty} k^2  \|X_k - X_k' \|_p$ and $\inf_m s_m^2 > 0$.
\end{thm}

The proof of Theorem~\ref{thmmdependent} is based on the following
decomposition. Let $Z_1,Z_2$ be independent unit Gaussian random
variables. Then
\begin{eqnarray*}
&& \sup_{x \in\mathbb{R}} \bigl|P (S_n/\sqrt{n} \leq x ) - \Phi
(x/s_{nm} ) \bigr|
\\
&&\qquad = \sup_{x \in\mathbb{R}} \bigl|P \bigl(S_{|m}^{(1)} \leq x
\sqrt{n} - S_{|m}^{(2)} \bigr) - P (Z_1
\overline{\sigma}_m \leq x - Z_2 {\overline{
\varsigma}}_m ) \bigr|
\\
&&\qquad\leq \mathbf{A} + \mathbf{B} + \mathbf{C},
\end{eqnarray*}
where $\mathbf{A},\mathbf{B},\mathbf{C}$ are defined as
\begin{eqnarray*}
\mathbf{A}& = &\sup_{x \in\mathbb{R}} \bigl|\mathbb{E} \bigl[P_{|\mathbb{F}_m}
\bigl(S_{|m}^{(1)}/\sqrt{n} \leq x - S_{|m}^{(2)}/
\sqrt{n} \bigr) - P_{|\mathbb{F}_m} \bigl(Z_1 \sigma_{|m}
\leq x - S_{|m}^{(2)}/\sqrt{n} \bigr) \bigr] \bigr|,
\\
\mathbf{B} &= &\sup_{x \in\mathbb{R}} \bigl|\mathbb{E} \bigl[P_{|\mathbb{F}_m}
\bigl(Z_1 \sigma_{|m} \leq x - S_{|m}^{(2)}/
\sqrt{n} \bigr) - P_{|\mathbb{F}
_m} \bigl(Z_1 \overline{\sigma}_{m} \leq x - S_{|m}^{(2)}/\sqrt {n} \bigr)
\bigr] \bigr|,
\\
\mathbf{C} &= & \sup_{x \in\mathbb{R}} \bigl|P \bigl(S_{|m}^{(2)}/
\sqrt{n} \leq x - Z_1 \overline{\sigma}_{m} \bigr) - P
(Z_2 {\overline {\varsigma}}_{m} \leq x - Z_1
\overline{\sigma}_{m} ) \bigr|.
\end{eqnarray*}
%
%
We will treat the three parts separately, and show that ${\mathbf A},
{\mathbf
B}, {\mathbf C} \leq\break \frac{c(\lambda,p)}{3}n^{-p/2 + 1}$, which proves
Theorem~\ref{thmmdependent}. As a brief overview, the proof consists
of the following steps:
\begin{longlist}[(a)]
\item[(a)] apply Esseen's smoothing inequality, and factor the
resulting characteristic function into a (conditional) product of
characteristic functions $\varphi_j(x)$ under the conditional
probability measure $P_{\mathbb{F}_m}$;
\item[(b)] use ideal metrics to control the distance between $\varphi
_j(x)$ and corresponding Gaussian versions under $P_{\mathbb{F}_m}$;
\item[(c)] based on Renyi's representation, control the (conditional)
characteristic functions $\varphi_j(x)$ under $P$;
\item[(d)] replace conditional variances under the overall probability
measure $P$.
\end{longlist}
One of the main difficulties arises from working under the conditional
measure~$P_{\mathbb{F}_m}$. For the proof, we require some additional
notation. In analogy to the filter $\theta_{k}^{(l, \prime)}$, we denote
with $\theta_{k}^{(l,*)}$,
%
\begin{eqnarray}
\label{defnstrichdepe2}
&& \theta_{k}^{(l,*)} = \bigl(
\epsilon_{k}, \epsilon_{k - 1},\ldots ,\epsilon_{k - l}',
\epsilon_{k - l - 1}',\epsilon_{k-l-2}',
\ldots \bigr).
\end{eqnarray}
We put $\theta_{k}^* = \theta_{k}^{(k,*)} =  (\epsilon_{k },
\epsilon_{k - 1},\ldots,\epsilon_{0}',\epsilon_{- 1}',\epsilon_{-
2}',\ldots )$ and ${X}_{k}^{(l,*)} = g (\theta_{k}^{(l,*)}
 )$, and in particular, we have ${X}_{k}^{*} = {X}_{k}^{(k,*)}$.
Similarly, let $ \{\epsilon_k'' \}_{k \in\mathbb{Z}}$ be
independent copies of $ \{\epsilon_k \}_{k \in\mathbb{Z}}$ and
$ \{\epsilon_k' \}_{k \in\mathbb{Z}}$. For $l \leq k$, we then
introduce the quantities $X_k^{(l,\prime\prime)}, X_k^{(l,**)}, X_k^{\prime\prime},
X_k^{(**)}$ in analogy to $X_k^{(l,\prime)}, X_k^{(l,*)}, X_k^{\prime}, X_k^{*}$.
This means that we replace every $\epsilon_k'$ with $\epsilon_k''$ at
all corresponding places. For $k \geq0$, we also introduce the $\sigma$-algebras
%
\begin{eqnarray}
 \mathcal{E}_k' &=& \sigma \bigl(\epsilon_j,
j \leq k \mbox{ and } j \neq0, \epsilon_0' \bigr)\quad\mbox{and}
\nonumber
\\[-8pt]
\\[-8pt]
\nonumber
\mathcal {E}_k^* &=& \sigma \bigl(\epsilon_j, 1
\leq j \leq k \mbox{ and } \epsilon_i', i \leq0 \bigr).
\end{eqnarray}
Similarly, we define $\mathcal{E}_k''$ and $\mathcal
{E}_k^{**}$.

Throughout the proofs, we make the following conventions:
\begin{longlist}[(4)]
\item[(1)] We do not distinguish between $N$ and $N + m'/2m$ since the
difference $m'/2m$ is not of any particular relevance for the proofs.
We use $N$ for both expressions.
\item[(2)] The abbreviations $I, \mathit{II}, \mathit{III}, \ldots,$ for expressions
(possible with some additional indices) vary from proof to proof.
\item[(3)] We use $\lesssim$, $\gtrsim$, ($\thicksim$) to denote
(two-sided) inequalities involving a multiplicative constant.
\item[(4)] If there is no confusion, we put $Y_j =
(2m)^{-1/2}Y_j^{(1)}$ for $j = 1,\ldots,N$ to lighten the notation,
particularly in part ${\mathbf A}$.
\item[(5)] We write [as in (\ref{defnsnm})] $\stackrel{\mathrm{def}}{=}$ if
we make definitions on the fly.
\end{longlist}

\subsubsection{Part $\mathbf{A}$}



The proof of part $\mathbf{A}$ is divided into four major steps. Some more
technical arguments are deferred to Sections~\ref{seczolo} and~\ref{secboundcondichar}.

\begin{pf}
For $L > 0$, put $\mathcal{B}_{L} =  \{{L}^{-1}\sum_{j = 1}^L
\sigma_{j|m}^2 \geq s_m^2/4 \}$, and denote with $\mathcal
{B}_{L}^{c}$ its complement. Since $S_{|m}^{(2)} \in\mathbb{F}_m$, we
obtain that
\begin{eqnarray*}
\mathbf{A} &= &\sup_{x \in\mathbb{R}} \bigl|\mathbb{E} \bigl[P_{\mathbb{F}_m}
\bigl(S_{|m}^{(1)}/\sqrt{n} \leq x - S_{|m}^{(2)}/
\sqrt{n} \bigr) - P_{\mathbb{F}_m} \bigl(Z_1 \sigma_{|m}
\leq x - S_{|m}^{(2)}/\sqrt{n} \bigr) \bigr]
\bigr|
\\
&\leq & \mathbb{E} \Bigl[\sup_{y \in
\mathbb{R}} \bigl|P_{\mathbb{F}_m}
\bigl(S_{|m}^{(1)}/\sqrt{n} \leq y \bigr) - P_{\mathbb{F}
_m}
(Z_1 \sigma_{|m} \leq y ) \bigr|\mathbh{1}(\mathcal
{B}_{N}) \Bigr] + 2P \bigl(\mathcal{B}_{N}^{c}
\bigr).
\end{eqnarray*}
Corollary~\ref{corsigreplaceissmall} yields that $P (\mathcal
{B}_{N}^{c} ) \lesssim n^{-p/2}N \lesssim n^{-p/2 + 1}$ since $N
\leq n$, and it thus suffices to treat
%
\begin{eqnarray}
\label{defnDeltam}
&& \Delta_{|m} \stackrel{\mathrm{def}} {=} \sup
_{y \in\mathbb{R}} \bigl|P_{\mathbb{F}
_m} \bigl(S_{|m}^{(1)}/
\sqrt{n} \leq y \bigr) - P_{\mathbb
{F}_m} (Z_1 \sigma_{|m}
\leq y ) \bigr|\mathbh{1}(\mathcal{B}_{N}).
\end{eqnarray}

\begin{longlist}[{}]
\item[\textit{Step} 1:] Berry--Esseen inequality.
Denote with $\Delta_{|m}^T$ the smoothed version of $\Delta_{|m}$
(cf. \cite{fellervolume2}) as in the classical approach. Since
$\sigma_{|m}^2 \geq s_m^2/4 > 0$ on the set $\mathcal{B}_N$ by
construction, the smoothing inequality (cf. \cite{fellervolume2}, Lemma~1, XVI.3) is applicable, and it thus suffices to treat
$\Delta_{|m}^T$. Let $\varphi_j(x) = \mathbb{E} [e^{\mathrm
{i}x Y_j} |
\mathbb{F}_m ]$, and put $T = n^{p/2-1} c_T$, where $c_T > 0$
will be
specified later. Due to the independence of $ \{Y_j \}_{1
\leq j \leq N}$ under $P_{|\mathbb{F}_m}$ and since $\mathbh
{1}(\mathcal
{B}_{N})\leq1$, it follows that
%
\begin{eqnarray}
\label{eqberrysmoothing}
&& \mathbb{E} \bigl[ \bigl|\Delta_{|m}^T \bigr| \bigr]
 \leq\int_{-T}^{T}\mathbb{E} \Biggl[ \Biggl|\prod
_{j = 1}^N \varphi_j (\xi/
\sqrt{N} ) - \prod_{j = 1}^N
e^{-\sigma_{j|m}^2 \xi^2/2N} \Biggr| \Biggr] \Big/|\xi| \,d \xi.
\end{eqnarray}
Put $t = \xi/\sqrt{N}$. Then $\prod_{j = 1}^N a_j - \prod_{j = 1}^N
b_j = \sum_{i = 1}^N  (\prod_{j = 1}^{i-1}b_j )(a_i -
b_i) (\prod_{j = i+1}^N a_j )$, where we use the convention
that $\prod_{j = 1}^{i-2}(\cdot) = \prod_{j = i + 2}^{N}(\cdot) = 1$ if $i-2
< 1$ or \mbox{$i + 2 > N$}. Hence we have
\begin{eqnarray*}
&& \prod_{j = 1}^N \varphi_j(t)
- \prod_{j = 1}^N e^{-\sigma_{j|m}^2
t^2/2}\\
&&\qquad = \sum
_{i = 1}^N \Biggl(\prod
_{j = 1}^{i-1}\varphi _j(t) \Biggr) \bigl(
\varphi_i(t) -e^{-\sigma_{i|m}^2 t^2/2} \bigr) \Biggl(\prod
_{j = i+1}^N e^{-\sigma_{j|m}^2 t^2/2} \Biggr).
\end{eqnarray*}
Note that both $ \{\varphi_j(t) \}_{1 \leq j \leq N}$ and
$ \{e^{-\sigma_{j|m}^2 t^2/2} \}_{1 \leq j \leq N}$ are
two-dependent sequences. Since $|\varphi_j(t)|$, $e^{-\sigma_{j|m}^2
t^2/2} \leq1$, it then follows by the triangle inequality,
stationarity and ``leave one out'' that
\begin{eqnarray*}
&& \Biggl\|\prod_{j = 1}^N \varphi_j(t)
- \prod_{j = 1}^N e^{-\sigma
_{j|m}^2 t^2/2}
\Biggr\|_1
\\
&&\qquad \leq\sum_{i = 1}^N \Biggl\|\prod
_{j =
1}^{i-2}e^{-\sigma_{j|m}^2 t^2/2} \Biggr\|_1 \bigl\|
\varphi_i(t) -e^{-\sigma_{i|m}^2 t^2/2} \bigl\|_1 \Biggl\|\prod
_{j = i+2}^N \bigl|\varphi_j(t) \bigr|
\Biggr\|_1
\\
&&\qquad \leq N \bigl\|\varphi_1(t) -e^{-\sigma_{1|m}^2 t^2/2} \bigr\|_1 \Biggl\|\prod
_{j = N/2}^{N-1} \bigl|\varphi_j(t) \bigr|
\Biggr\|_1
\\
&&\qquad\quad{}+ N \Biggl\|\prod_{j = 1}^{N/2-3}e^{-\sigma_{j|m}^2 t^2/2}
\Biggr\|_1 \bigl\|\varphi_1(t) -e^{-\sigma_{1|m}^2 t^2/2} \bigr\|_1
\\
&&\qquad\quad{}+ \Biggl\|\prod_{j = 1}^{N/2-3}e^{-\sigma_{j|m}^2 t^2/2}
\Biggr\|_1 \bigl\|\varphi_N(t) -e^{-\sigma_{N|m}^2 t^2/2} \bigr\|_1\\
&&\qquad =
I_N(\xi ) + \mathit{II}_N(\xi) + \mathit{III}_N(\xi).
\end{eqnarray*}

We proceed by obtaining upper bounds for $I_N(\xi), \mathit{II}_N(\xi)$ and
$\mathit{III}_N(\xi)$.

\item[\textit{Step} 2:]
Bounding $ \|\varphi_i(t) -e^{-\sigma
_{i|m}^2 t^2/2} \|_1$, $i \in\{1,N\}$. Let $Z_i$, $i \in\{1,N\}$
be two zero mean standard Gaussian random variables. Then
\begin{eqnarray*}
\bigl\|\varphi_i(t) -e^{-\sigma_{i|m}^2 t^2/2} \bigr\|_1 &\leq &  \bigl\|
\mathbb{E}_{\mathbb{F}_m} \bigl[\cos(tY_i) - \cos(t\sigma
_{i|m}Z_i ) \bigr] \bigr\|_1
\\
&&{}+ \bigl\|\mathbb{E}_{\mathbb{F}_m} \bigl[\sin(tY_i) - \sin(t\sigma
_{i|m}Z_i ) \bigr] \bigr\|_1.
\end{eqnarray*}
Due to the very nice analytical properties of $\sin(y), \cos(y)$, one
may reformulate the above in terms of ideal-metrics; cf. \cite{zolotarev1977} and Section~\ref{seczolo}. This indeed leads
to the desired bound
%
\begin{eqnarray}
\label{eqthmauxcharequationbound1}
&& \bigl\|\varphi_i(t) -e^{-\sigma_{i|m}^2 t^2/2} \bigr\|_1
\lesssim |t|^p m^{-p/2 + 1}.
\end{eqnarray}
The precise derivation is carried out in Section~\ref{seczolo} via
Lemmas~\ref{lemtrunccharcosandsin} and \ref{lemapproxchar},
and Corollary~\ref{corzolo}. Whether $i = 1$ or $i = N$ makes no difference.


\item[\textit{Step} 3:]
Bounding\vspace*{1.5pt} $ \|\prod_{j = N/2}^{N-1}
|\varphi_j(t) | \|_1$: in order to bound $ \|\prod_{j
= N/2}^{N-1}  |\varphi_j(t) | \|_1$, we require good
enough estimates for $ |\varphi_j(t) |$ where $0 \leq t <
1$. As already mentioned, we cannot directly follow the classical
approach. Instead, we use a refined version based on a conditioning
argument. To this end, let us first deal with $\varphi_j(t)$. Put
%
\begin{eqnarray}
\label{eqthmauxcharequationbound2.88}
&& \mathcal{G}_j^{(l)} = \sigma\big(\mathcal{E}_{(2j-2)m + l} \cup \{\epsilon_{(2j-1)m+1},\ldots,\epsilon_{2jm}\}\big)
\qquad \mbox{for $j,l \geq1$}.
\end{eqnarray}
We first consider the case $j = 1$. Introduce
%
\begin{eqnarray}
\mathit{IV}_{1}^{(l)}(m) &=& \sum
_{k = l+1}^m \bigl(X_k -
\mathbb{E}_{\mathbb
{F}_m} [X_k ] \bigr) + R_1\quad\mbox{and}
\nonumber
\\[-8pt]
\label{eqthmauxcharequationbound2.89}
\\[-8pt]
\nonumber
 V_{1}^{(l)}(m)  &=&  \mathit{IV}_1^{(l)}(m)
- \mathbb{E}_{\mathcal
{G}_1^{(l)}} \bigl[\mathit{IV}_1^{(l)}(m) \bigr].
\end{eqnarray}
Then
%
\begin{eqnarray}
\label{eqthmauxcharequationbound2.9}
&& \bigl|\varphi_1(t) \bigr| \leq\mathbb{E}_{\mathbb{F}_m} \bigl[
\bigl|\mathbb{E} _{\mathcal{G}_1^{(l)}} \bigl[e^{\mathrm
{i}t(2m)^{-1/2}V_1^{(l)}(m)} \bigr] \bigr| \bigr].
\end{eqnarray}
Clearly, this is also valid for $\varphi_j(t)$, $j = 2,\ldots,N$, with
corresponding $\mathcal{G}_j^{(l)}$ and $\mathit{IV}_j^{(l)}(m)$,
$V_j^{(l)}(m)$, defined analogously to (\ref
{eqthmauxcharequationbound2.89}). Let
%
\begin{eqnarray}
&& \varphi_j^{(l)}(x) = \mathbb{E}_{\mathcal{G}_j^{(l)}}
\bigl[e^{\mathrm{i}
x(m-l)^{-1/2}V_j^{(l)}(m)} \bigr],
\end{eqnarray}
and $\mathcal{J} =  \{j \dvtx  N/2 \leq j \leq N-1 \mbox{ and }2
\mbox{ divides }j\}$, and hence $\mathcal{J}$ denotes the set of all
even numbers between $N/2$ and $N-1$. Then
\begin{eqnarray*}
&& \prod_{j = N/2}^{N-1}\bigl|\varphi_j(t)
\bigr| \leq\prod_{j \in
\mathcal{J} }\mathbb{E}_{\mathbb{F}_m} \bigl[ \bigl|
\mathbb {E}_{\mathcal
{G}_j^{(l)}} \bigl[e^{\mathrm{i}t(2m)^{-1/2}V_j^{(l)}(m)} \bigr] \bigr| \bigr] = \prod
_{j \in\mathcal{J} }\mathbb{E}_{\mathbb{F}_m} \bigl[ \bigl|\varphi
_j^{(l)}(x) \bigr| \bigr],
\end{eqnarray*}
%
%
where $x = t\sqrt{(m-l)/2m}$. Note that $ \{V_j^{(l)}(m) \}
_{j \in\mathcal{J}}$ is a sequence\vspace*{1pt} of i.i.d. random variables,
particularly with respect to $P_{\mathbb{F}_m}$. Hence by independence and
Jensen's inequality, it follows from the above that
%
\begin{eqnarray}
 \Biggl\|\prod_{j = N/2}^{N-1} \bigl|
\varphi_j(t) \bigr| \Biggr\|_1  &\leq & \prod
_{j \in\mathcal{J} } \bigl\|\mathbb{E}_{\mathbb
{F}_m} \bigl[ \bigl|
\varphi_j^{(l)}(x) \bigr| \bigr] \bigr\|_1
\nonumber
\\[-8pt]
\label{eqthmauxcharequationbound2.905}
\\[-8pt]
\nonumber
 &\leq & \prod
_{j \in
\mathcal{J} } \bigl\|\varphi_j^{(l)}(x)
\bigr\|_1 = \Biggl\|\prod_{j \in\mathcal{J} } \bigl|
\varphi_j^{(l)}(x) \bigr| \Biggr\|_1.
\end{eqnarray}
We thus see that it suffices to deal with $\varphi_j^{(l)}(x)$. The
classical argument uses the estimate
\begin{eqnarray*}
&& \varphi \bigl(\xi/\sqrt{\sigma^2 n} \bigr) \leq e^{-5\xi^2/18n}
\qquad \mbox{for $\xi^2/n \leq c$, $c > 0$}
\end{eqnarray*}
for the characteristic function $\varphi$. Since in our case $\varphi
_j$ is random, we cannot use this estimate. Instead, we will use
Lemma~\ref{lemgeneralvarphibound}, which provides a similar result. In
order to apply it, set $J=  |\mathcal{J} |\geq N/8$,
%
\begin{eqnarray}
\label{eqthmauxcharequationbound2.91}
&& H_j = \frac{1}{\sqrt{m-l}}V_j^{(l)}(m)
\quad \mbox{and}\quad \mathcal{H}_j = \mathcal{G}_j^{(l)}.
\end{eqnarray}
For the applicability of Lemma~\ref{lemgeneralvarphibound}, we need
to verify that:
\begin{longlist}[(iii)]
\item[(i)] $\mathbb{E}_{\mathcal{H}_j} [H_j ] = 0$;
\item[(ii)] there exists a $u^- > 0$ such that $P (\mathbb
{E}_{\mathcal
{H}_j} [H_j^2 ] \leq u^- ) < 1/7$, uniformly for $j \in
\mathcal{J}$;
\item[(iii)] $ \|H_j \|_p \leq c_1$ uniformly for $j \in
\mathcal{J}$ and some $c_1 < \infty$.
\end{longlist}
Now (i) is true by construction. Claim (ii) is dealt with via
Lemma~\ref{lemauxHsigcomp}, which yields that
%
\begin{eqnarray}
\label{eqthmauxcharequationbound2.92}
&& P \bigl(\mathbb{E}_{\mathcal{H}_j} \bigl[H_j^2
\bigr] \leq\widehat {\sigma }_{m-l}^2 \bigr) \lesssim
\frac{1}{\sqrt{m-l}}.
\end{eqnarray}
%
%
Since $\widehat{\sigma}_{m-l}^2 \geq s_{m}^2/4$ for large enough
$m-l$ (say $m-l \geq K_0 > 0$) by Lemma~\ref
{lemsigexpressionsrelations}, we may set $0 < u^- = s_m^2/8 \leq
\widehat{\sigma}_{m-l}^2/2$. For showing (iii), it suffices to treat
the case $j = 1$. Note that (for $k\leq m$)
%
\begin{eqnarray}
\label{eqthmauxcharequationbound3}
&&
\quad\mathbb{E}_{\mathcal{G}_1^{(l)}} [X_k ] = \mathbb
{E}_{\mathcal{E}_l} \bigl[X_k - X_k^{(k-l,*)}
\bigr]\quad \mbox{and} \quad \mathbb{E}_{\mathbb
{F}_m} [X_k ] =
\mathbb{E}_{\mathbb{F}_m} \bigl[X_k-X_{k}^* \bigr].
\end{eqnarray}
By stationarity and the triangle and Jensen inequalities, 
we then have that
%
\begin{eqnarray}
\nonumber
\sqrt{m-l} \|H_j \|_p &\leq &  \Biggl\|
\sum_{k = l+1}^m X_k
\Biggr\|_p + \Biggl\|\sum_{k = l+1}^m
\mathbb{E}_{\mathcal
{E}_l} \bigl[X_k - X_{k}^{(k-l,*)}
\bigr] \Biggr\|_p
\nonumber
\\[-8pt]
\label{eqthmauxcharequationbound4}
\\[-8pt]
\nonumber
&&{}+ 2 \Biggl\|\sum_{k = l+1}^m \mathbb{E}_{\mathbb{F}_m}
\bigl[X_k - X_k^* \bigr] \Biggr\|_p + 2
\|R_1 \|_p.
\end{eqnarray}
Using Jensen's inequality and arguing similar to Lemma~\ref
{lemboundR1}, it follows that
\begin{eqnarray*}
\Biggl\|\sum_{k = l+1}^m\mathbb{E}_{\mathcal{G}_l}
\bigl[X_k - X_{k}^{(k-l,*)} \bigr] \Biggr\|_p
&\leq & \sum_{k = l+1}^m \bigl\|X_k -
X_k^{(k-l,*)} \bigr\|_p
\\
&\leq & \sum_{k = 1}^\infty k \bigl\|X_k -
X_k' \bigr\|_p < \infty.
\end{eqnarray*}
Similarly, using also Lemma~\ref{lemboundR1} to control $\|R_1\|_p$,
we obtain that
%
\begin{eqnarray}
\label{eqthmauxcharequationbound5}
&&\Biggl\|\sum_{k = l+1}^m
\mathbb{E}_{\mathbb{F}_m} [X_k ] \Biggr\|_p + \|R_1
\|_p < \infty.
\end{eqnarray}
By Lemma~\ref{lemwuoriginal}, we have $ \|\sum_{k = l+1}^m
X_k \|_p \lesssim\sqrt{m-l}$, and hence (iii) follows. We can
thus apply Lemma~\ref{lemgeneralvarphibound} with $u^- = s
_m^2/8$ and $J=  |\mathcal{J} | \geq N/8$, which yields
%
\begin{eqnarray}
\label{eqthmauxcharequationbound6}
&& \biggl\|\prod_{j \in\mathcal{J} } \bigl|\varphi_j^{(l)}(x)
\bigr| \biggr\|_1 \lesssim e^{-c_{\varphi,1}x^2 N/16} + e^{-\sqrt
{N/32}\log8/7}\qquad \mbox{for $x^2 < c_{\varphi,2}$,}
\end{eqnarray}
where $x = t\sqrt{(m-l)/2m}$. It is important to emphasize that both
$c_{\varphi,1}, c_{\varphi,2}$ do not depend on $l,m$ and are
strictly positive. Moreover, we find from (\ref
{eqthmauxcharequationbound2.92}) that $l$ can be chosen freely, as
long as $m-l$ is larger than $K_0$, which will be important in the next step.

\item[\textit{Step} 4:] Bounding\vspace*{1pt} and integrating $I_N(\xi),\mathit{II}_N(\xi
),\mathit{III}_N(\xi)$.

We first treat $I_N(\xi)$. Recall that $t = \xi/\sqrt{N}$, hence
\begin{eqnarray*}
&& |t|^p m^{-p/2 + 1} \lesssim|\xi|^p
n^{-p/2 + 1}N^{-1}.
\end{eqnarray*}
By (\ref{eqthmauxcharequationbound1}), (\ref
{eqthmauxcharequationbound2.905}) and (\ref
{eqthmauxcharequationbound6}), it then follows for $\xi^2 (m-l) <
c_{\varphi,2}n$ that
%
\begin{eqnarray}
&& I_N(\xi) \lesssim|\xi|^p n^{-p/2 + 1}
\bigl(e^{-c_{\varphi,1}\xi
^2 (m-l)/16m} + e^{-\sqrt{N/32}\log8/7} \bigr).
\end{eqnarray}
To make use of this bound, we need to appropriately select $l= l(\xi
)$. Recall that $N = n^{\lambda}$, $0 < \lambda\leq p/(2p + 2)$ by
assumption. Choosing
\begin{eqnarray*}
&& l(\xi) = \mathbh{1} \bigl(\xi^2 < n^{\lambda} c_{\varphi,2}
\bigr) + \biggl(m - \frac{c_{\varphi,2}n}{2\xi^2}\vee K_0 \biggr)\mathbh{1}
\bigl(\xi^2 \geq n^{\lambda} c_{\varphi,2} \bigr)
\end{eqnarray*}
and $c_T^2 < c_{\varphi,2}/K_0$, we obtain from the above that
%
\begin{eqnarray}
&& \int_{-T}^{T}I_N(\xi)/\xi \,d \xi
\lesssim n^{-p/2 + 1}.
\end{eqnarray}
In order to treat $\mathit{II}_N(\xi)$, let $N' = N/2 - 3$, and $\mathcal{B}_{N'} =  \{{N'}^{-1}\sum_{j = 1}^{N'} \sigma
_{j|m}^2 \geq s_m^2/4  \}$. Denote\vspace*{1pt} with $\mathcal{B}_{N'}^{c}$
its complement. Then by Corollary~\ref{corsigreplaceissmall}
(straightforward adaption is necessary) and (\ref
{eqthmauxcharequationbound1}), it follows that
%
\begin{eqnarray}
\nonumber
\mathit{II}_N(\xi) \mathbh{1} \bigl(|\xi|\leq N\bigr)&\leq &  N \bigl\|
\varphi_1(t) -e^{-\sigma_{1|m}^2 \xi^2/2} \bigr\|_1 \Biggl\|\prod
_{j =
1}^{N'}e^{-\sigma_{j|m}^2 \xi^2/2}\mathbh{1} (\mathcal
{B}_{N'} ) \Biggr\|_1
\\
&&\label{eqthmauxcharequationbound9}
{}+ N \bigl\|\varphi_1(t) -e^{-\sigma
_{1|m}^2 \xi^2/2} \bigr\|_1 P
\bigl(\mathcal{B}_{N'}^c \bigr)
\\
\nonumber
&\lesssim & |\xi|^p n^{-p/2 + 1} e^{-s_{m}^2 \xi^2/16} + |
\xi|^p n^{-p/2 + 1} N n^{-p/2}.
\end{eqnarray}
Similarly, using $ \|\varphi_1(t) -e^{-\sigma_{1|m}^2 \xi
^2/2} \|_1 \leq2$ one obtains
%
\begin{eqnarray}
\label{eqthmauxcharequationbound10}
&& \mathit{II}_N(\xi) \mathbh{1} \bigl(|\xi| > N\bigr) \lesssim|
\xi|^p n^{-p/2 + 1} e^{-s_{m}^2 \xi^2 N/16} + n^{-p/2}
N^{2}.
\end{eqnarray}
%
Hence employing (\ref{eqthmauxcharequationbound9}) and (\ref
{eqthmauxcharequationbound10}) yields
%
\begin{eqnarray}
\nonumber
\int_{-T}^{T}\mathit{II}_N(\xi)/\xi
\,d\xi &\lesssim &  n^{-p/2 + 1} \int_{|\xi
| \leq N} |
\xi|^{p-1} \bigl(e^{-s_{m}^2 \xi^2/16} + n^{-p/2} N \bigr) \,d\xi
\\
&&{}+ \int_{N < |\xi| \leq T} \bigl(n^{-p/2 + 1} |
\xi|^{p-1} e^{-s_{m}^2 \xi^2/16} + n^{-p/2} N^2
\xi^{-1} \bigr) \,d\xi
\nonumber
\\[-8pt]
\\[-8pt]
\nonumber
& \lesssim &  n^{-p/2 + 1} + n^{-p/2+1} n^{-p/2} N^{p+1}
+ n^{-p/2} N^2 \log T \\
\nonumber
&\lesssim &  n^{-p/2 + 1},
\end{eqnarray}
since $N = n^{\lambda}$, $0 < \lambda\leq p/(2p + 2)$ by assumption.
Similarly, one obtains the same bound for $\mathit{III}_N(\xi)$. This completes
the proof of part ${\mathbf A}$.\quad\qed
\end{longlist}
\noqed\end{pf}

\subsubsection{Part $\mathbf{B}$}
\mbox{}
\begin{pf}
Let
\begin{eqnarray*}
&& \Delta^{(2)}(x) \stackrel{\mathrm{def}} {=} \mathbb{E} \bigl[P_{|\mathbb
{F}_m}
\bigl(Z_1 \sigma_{|m} \leq x - S_{|m}^{(2)}/
\sqrt{n} \bigr) - P_{|\mathbb{F}
_m} \bigl(Z_1 \overline{
\sigma}_{m}\leq x - S_{|m}^{(2)}/\sqrt {n} \bigr)
\bigr].
\end{eqnarray*}
Recall that $\mathcal{B}_{N} =  \{{N}^{-1}\sum_{j = 1}^N \sigma
_{j|m}^2 \geq s_m^2/4  \}$ and $P (\mathcal{B}_N^c )
\lesssim n^{-p/2 + 1}$ by Corollary~\ref{corsigreplaceissmall}.
Using properties of the Gaussian distribution, it follows that
\begin{eqnarray*}
{\mathbf B} &\leq & \sup_{x\in\mathbb{R}} \bigl|\mathbb{E} \bigl[
\Delta^{(2)}(x) \mathbh{1} (\mathcal{B}_{N}) \bigr] \bigr| + \sup
_{x\in\mathbb{R}} \bigl|\mathbb{E} \bigl[\Delta^{(2)}(x) \mathbh{1}\bigl(
\mathcal{B}_{N}^c\bigr) \bigr] \bigr|
\\
&\lesssim & \mathbb{E} \bigl[ \bigl|1/\sigma_{|m} - 1/\overline{
\sigma}_{m} \bigr|\mathbh{1} (\mathcal{B}_N) \bigr] +
n^{-p/2 + 1}.
\end{eqnarray*}
Using $(a - b)(a+b) = a^2 - b^2$, H\"{o}lders inequality and Lemma~\ref
{lemsigexpansion}, we obtain that
\begin{eqnarray*}
&& \mathbb{E} \bigl[ \bigl|1/\sigma_{|m} - 1/\overline{\sigma
}_{m} \bigr|\mathbh{1} (\mathcal{B}_{N}) \bigr] \lesssim \bigl\|
\sigma_{|m}^2 - \overline {\sigma}_{m}^2
\bigr\|_{p/2} \lesssim n^{-p/2 + 1}.
\end{eqnarray*}
Hence we conclude that ${\mathbf B} \lesssim n^{-p/2 + 1}$.
\end{pf}

\subsubsection{Part $\mathbf{C}$}
\mbox{}
\begin{pf}
Due to the independence of $Z_1, Z_2$, we may rewrite ${\mathbf C}$ as
\begin{eqnarray*}
&& {\mathbf C} = \sup_{x \in\mathbb{R}} \bigl|\Phi \bigl( \bigl(x -
S_{|m}^{(2)}/\sqrt{n} \bigr)/\overline{\sigma}_{m}
\bigr) -\Phi \bigl( (x - Z_2 \overline{\varsigma}_{m} )/
\overline {\sigma}_{m} \bigr) \bigr|,
\end{eqnarray*}
where $\Phi(\cdot)$ denotes the c.d.f. of a standard normal
distribution. This induces a ``natural'' smoothing. The claim now
follows by repeating the same arguments as in part ${\mathbf A}$. Note
however, that the present situation is much easier to handle, due to
the already smoothed version, and since $Y_k^{(2)}$, $k = 1,\ldots,N$ is a
sequence of independent random variables. Alternatively, one may also
directly appeal to the results in \cite{dedeckerminimaldistance2009}.
\end{pf}

\subsection{Proof of Theorem~\texorpdfstring{\protect\ref{thmberry}}{2.2}}\label
{secproofofmaintheorem}

The proof of Theorem~\ref{thmberry} mainly consists of constructing a
good $m$-dependent approximation and then verifying the conditions of
Theorem~\ref{thmmdependent}. To this end, set $m = c n^{3/4}$ for
some $c > 0$, and note that $1/4 < p/(2p + 2)$ for $p \in(2,3]$. Let
$\mathcal{E}_k^{m} = \sigma (\epsilon_j, k-m +1 \leq j \leq
k )$, and define the approximating sequence as
%
\begin{eqnarray}
X_k^{(\leq m)} &=& \mathbb{E} \bigl[X_k |
\mathcal{E}_k^{m} \bigr]\quad \mbox{and}
\nonumber
\\[-8pt]
\\[-8pt]
\nonumber
X_k^{(>m)} &=& X_k - X_k^{(\leq m)}
= X_k - \mathbb{E} \bigl[X_k | \mathcal{E}_k^{m}
\bigr].
\end{eqnarray}
We also introduce the corresponding partial sums as
%
\begin{eqnarray}
&& S_n^{(\leq m)} = \sum_{k = 1}^n
X_k^{(\leq m)},\qquad  S_n^{(> m)} = \sum
_{k = 1}^n X_k^{(> m)}.
\end{eqnarray}
Further,\vspace*{1pt} let $s_n^2 = n^{-1} \|S_n \|_2^2$ and $s_{nm}^2 =
n^{-1} \|S_n^{(\leq m)} \|_2 = \overline{\sigma}_m^2 +
\overline{\varsigma}_m^2$. We require the following auxiliary result
(Lemma~5.1 in \cite{hormann2009}).

\begin{lem}\label{lemsiegidecomposition}
For every $\delta> 0$, every $m,n \geq1$ and every $x \in\mathbb
{R}$, the
following estimate holds:
\begin{eqnarray*}
&& \bigl |P (S_n /\sqrt{n}\leq x s_n ) - \Phi(x) \bigr| \\
&&\qquad\leq
A_0(x,\delta) + A_1(m,n,\delta)
\\
&&\qquad\quad{}+ \max \bigl\{A_2(m,n,x,\delta) + A_3(m,n,\delta),
A_4(m,n,x,\delta) + A_5(m,n,x,\delta) \bigr\},
\end{eqnarray*}
where:
\begin{eqnarray*}
 A_0(x,\delta) &=& \bigl|\Phi(x) - \Phi(x + \delta) \bigr|;
\\
 A_1(m,n,\delta) &=&  P \bigl(\bigl|S_n - S_n^{(\leq m)}\bigr|
\geq\delta s_n \sqrt{n} \bigr);
\\
 A_2(m,n,x,\delta) &=& \bigl|P \bigl(S_n^{(\leq m)} \leq(x
+ \delta )s_n\sqrt{n} \bigr) - \Phi \bigl((x+\delta)s_n/s_{nm}
\bigr) \bigr|;
\\
 A_3(m,n,x,\delta) &=& \bigl|\Phi \bigl((x+\delta)s_n/s_{nm}
\bigr) - \Phi(x+\delta) \bigr|;
\\
 A_4(m,n,x,\delta) &=& A_2(m,n,x,-\delta)\quad \mbox{and}\quad
A_5(m,n,x,\delta) = A_3(m,n,x,-\delta).
\end{eqnarray*}
\end{lem}

\begin{pf*}{Proof of Theorem~\ref{thmberry}}
As a preparatory result, note that
%
\begin{eqnarray}
\label{eqthmgenvarid}&& n s_n^2 = n s^2 + \sum
_{k \in\mathbb{Z}}\bigl(n \wedge|k|\bigr)\mathbb {E} [X_0
X_k ].
\end{eqnarray}
Using the same arguments as in Lemma~\ref
{lemsigexpressionsrelations}, it follows that $n s_n = n s^2 +
\mathcal{O} (1 ) > 0$. By the properties of Gaussian distribution,
\begin{eqnarray*}
&& \sup_{x \in\mathbb{R}} \bigl|\Phi \bigl(x/\sqrt{s^2} \bigr) - \Phi
(x/\sqrt{s_n} ) \bigr| \lesssim n^{-1},
\end{eqnarray*}
and we may thus safely interchange $s_n^2$ and $s^2$. We first deal
with $A_1(m,n,\delta)$. For $j \in\mathbb{Z}$, denote with $\mathcal
{P}_j
(X_k^{(>m)} )$ the projection operator
%
\begin{eqnarray}
&& \mathcal{P}_j \bigl(X_k^{(>m)} \bigr) =
\mathbb{E} \bigl[X_k^{(>m)} |\mathcal {E}_j
\bigr] - \mathbb{E} \bigl[X_k^{(>m)} |\mathcal
{E}_{j-1} \bigr].
\end{eqnarray}
Proceeding as in the proof of Lemma~3.1 in \cite{jirakdarling2013},
it follows that for $k \geq0$,
%
\begin{eqnarray}
\label{eqthmgen4}
&&
\bigl\|\mathcal{P}_0 \bigl(X_k^{(>m)}
\bigr) \bigr\|_p \leq2 \min \Biggl\{ \bigl\|X_k -
X_k' \bigr\|_p, \sum
_{l = m}^{\infty} \bigl\|X_l -
X_l' \bigr\|_p \Biggr\}.
\end{eqnarray}
An application of Theorem~1 in \cite{sipwu} now yields that
%
\begin{eqnarray}
&& n^{-1/2} \bigl\|S_n^{(>m)} \bigr\|_p \leq c(p)
\sum_{k = 1}^{\infty} \bigl\|\mathcal{P}_0
\bigl(X_k^{(>m)} \bigr) \bigr\|_p
\end{eqnarray}
for some absolute constant $c(p)$ that only depends on $p$. By (\ref
{eqthmgen4}), it follows that the above is of magnitude
%
\begin{eqnarray}
\label{eqthmgen6}
&& \sum_{k = L}^{\infty}
L^{-2} k^2 \bigl\|X_k - X_k'
\bigr\|_p + L \sum_{k = m}^{\infty}
m^{-2} k^2 \bigl\|X_k - X_k'
\bigr\|_p \lesssim L^{-2} + Lm^{-2}.\hspace*{-16pt}
\end{eqnarray}
Setting $L = m^{2/3}$, we obtain the bound $\mathcal{O}
(m^{-4/3} )
= \mathcal{O} (n^{-1} )$. We thus conclude from the Markov
inequality that
\begin{eqnarray*}
&& P \bigl(\bigl|S_n - S_n^{(\leq m)}\bigr| \geq\delta
s_n \sqrt{n} \bigr) = P \bigl(\bigl|S_n^{(> m)}\bigr| \geq
\delta s_n \sqrt{n} \bigr) \lesssim(\delta n)^{-p},
\end{eqnarray*}
hence
%
\begin{eqnarray}
\label{eqthmgena1}
&& A_1(m,n,\delta) \lesssim(\delta n)^{-p}.
\end{eqnarray}
Note that a much sharper bound can be obtained via moderate deviation
arguments (cf. \cite{grama1997}), but the current one is sufficient
for our needs, and its deviation requires fewer computations. Next, we
deal with $A_2(m,n,x,\delta)$. The aim is to apply Theorem~\ref
{thmmdependent} to obtain the result. In order to do so, we need to
verify Assumption~\ref{assdependence}(i)--(iii) for $X_k^{(\leq
m)}$.
\begin{longlist}[{}]
\item[\textit{Case} (i):] Note first that $\mathbb{E}
[X_k^{(\leq m)} ] = \mathbb{E}
 [X_k ] = 0$. Moreover, Jensen's inequality gives
\begin{eqnarray*}
&& \bigl\|X_k^{(\leq m)} \bigr\|_p = \bigl\|\mathbb{E}
\bigl[X_k | \mathcal{E}_k^m \bigr]
\bigr\|_p \leq \|X_k \|_p < \infty.
\end{eqnarray*}
Hence Assumption~\ref{assdependence}(i) is valid.
\item[\textit{Case} (ii):] Note that we may assume $k \leq m$, since otherwise
$ (X_k^{(\leq m)} )' - X_k^{(\leq m)} = 0$, and
Assumption~\ref{assdependence}(ii) is trivially true. Put
\begin{eqnarray*}
&&\mathcal{E}_k^{(m,\prime)} = \sigma \bigl(\epsilon_j,
k-m +1 \leq j \leq k, j \neq0, \epsilon_0' \bigr).
\end{eqnarray*}
Since $\mathbb{E} [X_k |\mathcal{E}_k^{m} ]' = \mathbb
{E}
[X_k' |\mathcal{E}_k^{(m,\prime)} ]$, it follows that
%
\begin{eqnarray}
\nonumber
\bigl(X_k^{(\leq m)} \bigr)' -
X_k^{(\leq m)} &=&  \mathbb{E}_{\mathcal
{E}_k^{(m,\prime)}}
\bigl[X_k' \bigr] - \mathbb{E}_{\mathcal
{E}_k^{m}}
[X_k ]
\\
&=& \mathbb{E}_{\mathcal{E}_k^{(m,\prime)}} \bigl[X_k' -
X_k \bigr] + \mathbb{E}_{\mathcal{E}_k^{(m,\prime)}} [X_k ] - \mathbb
{E}_{\mathcal
{E}_k^{m}} [X_k ]
\\
&=&
\nonumber
\mathbb{E}_{\mathcal
{E}_k^{(m,\prime)}} \bigl[X_k' -
X_k \bigr] + \mathbb{E}_{\mathcal
{E}_k^{m}} \bigl[X_k'
\bigr] - \mathbb{E}_{\mathcal{E}_k^{m}} [X_k ]
\\
& =& \mathbb{E} _{\mathcal{E}_k^{(m,\prime)}} \bigl[X_k' -
X_k \bigr] + \mathbb {E}_{\mathcal
{E}_k^{m}} \bigl[X_k'
- X_k \bigr].
\end{eqnarray}
Hence by Jensen's inequality $ \| (X_k^{(\leq m)} )' -
X_k^{(\leq m)} \|_p \leq2  \|X_k - X_k' \|_p$, which
gives the claim.
\item[\textit{Case} (iii):] We have $X_k^{(\leq m)} = \mathbb
{E} [X_k^{(m,*)}
|\mathcal{E}_k ]$. Then
%
\begin{eqnarray}
\nonumber
\bigl\|X_k^{(> m)} \bigr\|_p &=& \bigl\|
\mathbb{E} \bigl[X_k - X_k^{(m,*)} \bigr|
\mathcal{E}_k \bigr] \bigr\|_p \leq \bigl\|X_k -
X_k^{(m,*)} \bigr\|_p
\nonumber
\\[-8pt]
\label{eqthmgen7}
\\[-8pt]
\nonumber
&\leq &  m^{-2}\sum_{l = m}^{\infty}l^2
\bigl\| X_l - X_l' \bigr\|_p \lesssim
m^{-2}.
\end{eqnarray}
By the Cauchy--Schwarz, triangle and Jensen inequalities, we have
\begin{eqnarray*}
&& \bigl|\mathbb{E} [X_k X_0 ] - \mathbb{E}
\bigl[X_k^{(\leq
m)} X_0^{(\leq
m)} \bigr] \bigr|
\\
&&\qquad\leq \|X_0 \|_2 \bigl\|X_k^{(>
m)}
\bigr\|_2 + \|X_k \|_2 \bigl\|X_0^{(> m)}
\bigr\|_2 + \bigl\|X_0^{(> m)} \bigr\|_2
\bigl\|X_k^{(> m)} \bigr\|_2.
\end{eqnarray*}
By (\ref{eqthmgen7}), this is of the magnitude $\mathcal{O}
(m^{-2} )$. We thus conclude that
%
\begin{eqnarray}
&& \Biggl|\sum_{k = 0}^m\mathbb{E} [X_k
X_0 ] - \sum_{k =
0}^m
\mathbb{E} \bigl[X_k^{(\leq m)} X_0^{(\leq m)}
\bigr] \Biggr| \lesssim m^{-1}.
\end{eqnarray}

On the other hand, we have
\begin{eqnarray*}
&& \biggl|\sum_{k > m} \mathbb{E} [X_k
X_0 ] \biggr| \leq\sum_{k > m}
\|X_0 \|_2 \bigl\|X_k^* - X_k
\bigr\|_2 \leq\frac
{1}{m}\sum_{k > m}k^2
\bigl\|X_k - X_k' \bigr\|_2
\|X_0 \|_2 \lesssim\frac{1}{m}.
\end{eqnarray*}
This yields
%
\begin{eqnarray}
\label{eqthmgenvarest}
&& \Biggl|\sum_{k \in\mathbb{Z}}^m\mathbb{E}
[X_k X_0 ] - s^2 \Biggr| \lesssim
\frac{1}{m},
\end{eqnarray}
which gives (iii) for large enough $m$. Since $c n^{3/4}$, we see
that we may apply Theorem~\ref{thmmdependent} which yields
%
\begin{eqnarray}
\label{eqthmgena2}
&& \sup_{x \in\mathbb{R}} A_2(m,n,x,\delta)
\lesssim n^{-p/2 + 1}.
\end{eqnarray}
Next, we deal with $A_3(m,n,x,\delta)$. Properties of the Gaussian
distribution function give
\begin{eqnarray*}
&& \sup_{x \in\mathbb{R}}A_3(m,n,x,\delta) \lesssim\delta+
\bigl|s_n^2 - s_{nm}^2 \bigr|.
\end{eqnarray*}
However, by the Cauchy--Schwarz inequality and (\ref{eqthmgen6}),
it follows that
%
\begin{eqnarray}
&& \bigl|s_n^2 - s_{nm}^2 \bigr| \leq
n^{-1} \bigl\|S_n^{(> m)} \bigr\|_2
\bigl\|S_n + S_n^{(\leq m)} \bigr\|_2 \lesssim
m^{-4/3} \lesssim n^{-1},
\end{eqnarray}
and we thus conclude that
%
\begin{eqnarray}
\label{eqthmgena3}
&& \sup_{x \in\mathbb{R}}A_3(m,n,x,\delta)
\lesssim\delta+ n^{-1}.
\end{eqnarray}
Finally, setting $\delta= n^{-1/2}$, standard arguments involving the
Gaussian distribution function yield that
%
\begin{eqnarray}
\label{eqthmgena0}
&& \sup_{x \in\mathbb{R}}A_0(x,\delta) \lesssim
\delta= n^{-1/2}.
\end{eqnarray}
Piecing together (\ref{eqthmgena1}), (\ref{eqthmgena2}),
(\ref{eqthmgena3}) and (\ref{eqthmgena0}), Lemma~\ref
{lemsiegidecomposition} yields
%
\begin{eqnarray}
\label{eqthmgena10}
&& \sup_{x \in\mathbb{R}} \bigl|P (S_n /\sqrt{n}\leq x
) - \Phi (x/s_n ) \bigr|\lesssim n^{-p/2 + 1},
\end{eqnarray}
which completes the proof.\quad\qed
\end{longlist}
\noqed\end{pf*}

\subsection{Proof of Theorem~\texorpdfstring{\protect\ref{thmLp}}{2.3}}\label{secprooflp}
Recall that
\begin{eqnarray*}
&& \Delta_n(x) = \Bigl|P \Bigl(S_n \leq x \sqrt{n
s_n^2} \Bigr) - \Phi (x) \Bigr| \quad\mbox{and}\quad
\Delta_n = \sup_{x \in\mathbb{R}} \Delta_n(x).
\end{eqnarray*}
We first consider the case $q > 1$. Using Theorem~\ref{thmberry}, we have
%
\begin{eqnarray}
\label{eqthmLPboundrio2}
&& \int_{\mathbb{R}} \bigl|\Delta_n(x)
\bigr|^q \,dx \leq\Delta _n^{q-1}\int
_{\mathbb{R}} \bigl|\Delta_n(x) \bigr| \,dx \lesssim
n^{-({q-1})/{2}} \int_{\mathbb{R}} \bigl|\Delta_n(x) \bigr|\, dx.
\end{eqnarray}
In order to bound $\int_{\mathbb{R}} |\Delta_n(x) |\,dx$, we apply
\cite{dedeckerriomean2008}, Theorem~3.2, which will give us the bound
%
\begin{eqnarray}
\label{eqthmLPboundrio3}
&& \int_{\mathbb{R}} \bigl|\Delta_n(x) \bigr| \,dx
\lesssim\frac
{1}{\sqrt{n}}.
\end{eqnarray}
To this end, we need to verify that
%
\begin{eqnarray}
\nonumber
&&\sum_{k > 0} \Biggl(
\bigl\|X_0^2 \vee1 \bigl(\mathbb{E} \bigl[X_k^2
- \mathbb{E} \bigl[X_k^2\bigr] |\mathcal{E}_0
\bigr] \bigr) \bigr\|_1 \\
\label{eqRioL1}
&&\qquad{}+ \frac{1}{k}\sum
_{i = 1}^k \bigl\|X_{-i} X_0
\mathbb{E} \bigl[X_k^2 - \mathbb {E}\bigl[X_k^2
\bigr] | \mathcal{E}_0 \bigr] \bigr\|_1 \Biggr)<\infty\quad\mbox{and}
\\
\nonumber
&& \sum_{k > 0} \frac{1}{k} \sum
_{i = \lfloor k/2
\rfloor}^k \bigl\||X_0| \vee1
\mathbb{E} \bigl[X_i X_k^2 - \mathbb{E}
\bigl[X_i X_k^2\bigr] |
\mathcal{E}_0 \bigr] \bigr\|_1 < \infty.
\end{eqnarray}
Applying the H\"{o}lder, Jensen and triangle inequalities, we get
\begin{eqnarray*}
\bigl\|X_0^2 \vee1 \bigl(\mathbb{E} \bigl[X_k^2
- \mathbb {E}\bigl[X_k^2\bigr] |\mathcal
{E}_0 \bigr] \bigr) \bigr\|_1 &\leq &  \bigl\|X_0^2
\vee1 \bigr\|_2 \bigl\|X_k - X_k^* \bigr\|_4
\bigl\|X_k + X_k^* \bigr\|_4
\\
& \lesssim  & k \bigl\|X_k - X_k'
\bigr\|_4.
\end{eqnarray*}
Similarly, with $\mathcal{E}_{-i} = \sigma (\epsilon_k, k
\leq-i  )$, we obtain that
\begin{eqnarray*}
\bigl\|X_{-i} X_0 \mathbb{E} \bigl[X_k^2
- \mathbb{E}\bigl[X_k^2\bigr] | \mathcal
{E}_0 \bigr] \bigr\|_1 &\lesssim &  \|X_{-i}
\|_4 \bigl\| \mathbb{E} [X_{0}|\mathcal{E}_{-i}]
\bigr\|_4 k \bigl\|X_k - X_k'
\bigr\|_4
\\
&\lesssim &  \bigl\|X_i - X_i^* \bigr\|_4 k
\bigl\|X_k - X_k' \bigr\|_4
\\
&\lesssim &  i \bigl\|X_i - X_i'
\bigr\|_4 k \bigl\|X_k - X_k'
\bigr\|_4.
\end{eqnarray*}
In the same manner, we get that
\begin{eqnarray*}
\bigl\||X_0| \vee1 \mathbb{E} \bigl[X_i
X_k^2 - \mathbb{E}\bigl[X_i
X_k^2\bigr] | \mathcal{E}_0 \bigr]
\bigr\|_1 &\lesssim &  \bigl\|X_i - X_i^* \bigr\|
_4 + \bigl\|X_k - X_k^* \bigr\|_4
\\
&\lesssim &  i \bigl\|X_i - X_i'
\bigr\|_4 + k \bigl\|X_k - X_k'
\bigr\|_4.
\end{eqnarray*}
Combining all three bounds, the validity of (\ref{eqRioL1}) follows,
and hence (\ref{eqthmLPboundrio3}). For~(\ref{eqthmLPboundrio2}), we
thus obtain
\begin{eqnarray*}
&& \int_{\mathbb{R}} \bigl|\Delta_n(x) \bigr|^q\, dx
\lesssim n^{-{(q-1)}/{2} - 1/2} \lesssim n^{-q/2},
\end{eqnarray*}
which completes the proof for $q >1$. For $q = 1$, we may directly
refer to \cite{dedeckerriomean2008}, Theorem~3.2, using the above bounds.

\subsection{Proof of Theorem~\texorpdfstring{\protect\ref{thmnonunif}}{2.4}}\label
{secthmproofnonunif}

For the proof, we require the following result; cf.~\cite{petrovbook1995}, Lemma~5.4.

\begin{lem}\label{lemosipov}
Let $Y$ be a real-valued random variable. Put
\begin{eqnarray*}
&& \widetilde{\Delta} = \sup_{x \in \mathbb{R}}\bigl|P(Y \leq x ) - \Phi (x ) \bigr|,
\end{eqnarray*}
and assume that $\|Y\|_q < \infty$ for $q > 0$ and $0 \leq \widetilde{\Delta} \leq e^{-1/2}$. Then
\begin{eqnarray*}
&& \bigl|P (Y \leq x ) - \Phi (x ) \bigr| \leq
\frac{c(q)\widetilde{\Delta}  ( \log1/\widetilde{\Delta})^{q/2} + \lambda_q}{1 + |x|^{q}}
\end{eqnarray*}
for all $x$, where $c(q)$ is a positive constant depending only on $q$, and
\begin{eqnarray*}
&& \lambda_q = \biggl|\int_{\mathbb{R}} |x|^q \,d
\Phi(x) - \mathbb{E}\bigl[|Y|^q\bigr] \biggr|.
\end{eqnarray*}
\end{lem}

Consider first the case where $|x| \leq c_0 \sqrt{\log n}$, for $ c_0 > 0$ large enough (see below). Then by the Markov inequality and Lemma \ref{lemwuoriginal}, it follows that
\begin{eqnarray}
\label{eq_truncated_be}
&& \Bigl|P\Bigl(S_n \mathbh{1}\bigl(|S_n|\leq n\bigr) \leq x\sqrt{n s_n^2}\Bigr) - P\Bigl(S_n \leq x\sqrt{n s_n^2}\Bigr)\Bigr|\lesssim n^{-p/2}.
\end{eqnarray}
Combining Theorem~\ref{thmberry} with (\ref{eq_truncated_be}) and Lemma \ref{lemosipov}, we see that it suffices to consider $\lambda_p$ with $Y = S_n \mathbh{1}(|S_n|\leq n)$. Using again Theorem~\ref{thmberry} together with  (\ref{eq_truncated_be}), standard tail bounds for the Gaussian distribution and elementary computations give
%
\begin{eqnarray}
\label{eqthmnonunif2}
\quad\lambda_p & \lesssim &  n^{-(p\wedge3)/2 + 1} (\log n)^{p/2} + \int_{c_0 \sqrt{\log n}}^{n} x^{p-1}P \Bigl(|S_n| \geq x
\sqrt{n s_n^2} \Bigr)\, dx.
\end{eqnarray}
According to a Fuk--Nagaev-type inequality for dependent sequences in \cite{Wufuknagaev}, Theorem~2,
if it holds that
%
\begin{eqnarray}
\label{eqthmnonunif3}
&& \sum_{k = 1}^{\infty}
\bigl(k^{p/2 -1} \bigl\|X_k - X_k'
\bigr\|_p^p \bigr)^{{1}/({p+1})} < \infty,
\end{eqnarray}
then for large enough $c_0 > 0$ and $x \geq c_0 \sqrt{\log n}$ we get
%
\begin{eqnarray}
\label{eqthmnonunif4}
&& P\Bigl(|S_n| \geq x\sqrt{n s_n^2} \Bigr) \lesssim n^{-p/2 + 1} x^{-p},
\end{eqnarray}
and hence,
\begin{eqnarray}
\label{eqthmnonunif5}
&& \int_{c_0 \sqrt{\log n}}^{n} x^{p-1}P \Bigl(|S_n| \geq x
\sqrt{n s_n^2} \Bigr)\, dx \lesssim n^{-p/2 + 1} (\log n)^{p/2} \log n.
\end{eqnarray}
However, setting $a_k = k^{-{1}/{2} - {1}/({3(p+1)})}$, an application of the Cauchy--Schwarz inequality yields\vspace*{-3pt}
\begin{eqnarray*}
&& \Biggl(\sum_{k = 1}^{\infty}
\bigl(k^{p/2 -1} \bigl\|X_k - X_k'
\bigr\|_p^p \bigr)^{{1}/({p+1})} \Biggr)^2\\[-2pt]
 &&\qquad\leq
\sum_{k = 1}^{\infty} a_k^2
\sum_{k = 1}^{\infty} a_k^{-2}
k^{({p - 2})/({p+1})} \bigl\|X_k - X_k'\bigr\|
_p^{{2p}/({p+1})}
\\[-2pt]
&&\qquad\lesssim  \sum_{k = 1}^{\infty} k^{2}
\bigl\|X_k - X_k'\bigr\|_p < \infty
\end{eqnarray*}
by Assumption~\ref{assmaindependence}. Hence (\ref{eqthmnonunif3})
holds, and thus (\ref{eqthmnonunif4}) and (\ref{eqthmnonunif5}). To complete the proof, it remains to treat the case $|x| > c_0 \sqrt{\log n}$.
But in this case, we may directly appeal to (\ref{eqthmnonunif4}) which gives the result.

\subsection{Proof of main lemmas}\label{secmainlemmas}

\subsubsection{Bounding conditional characteristic functions and
variances}\label{secboundcondichar}

Suppose we have a sequence of random variables $ \{H_j \}_{1
\leq j \leq J}$ and a sequence of filtrations $ \{\mathcal
{H}_j \}_{1 \leq j \leq J}$, such that both $ \{\mathbb{E}
_{\mathcal{H}_j} [H_j^2 ] \}_{1 \leq j \leq J}$ and
$ \{\mathbb{E}_{\mathcal{H}_j} [|H_j|^p ] \}_{1
\leq j
\leq J}$ are independent sequences. Note that this does not
necessarily mean that $ \{H_j \}_{1 \leq j \leq J}$ is
independent, and indeed this is not the case when we apply Lemma~\ref{lemgeneralvarphibound} in step~4 of the proof of part ${\mathbf A}$.
Introduce the conditional characteristic function
%
\begin{eqnarray}
&& \varphi_j^{\mathcal{H}}(x) = \mathbb{E} \bigl[\exp(\mathrm{i}x
H_j) | \mathcal {H}_j \bigr].
\end{eqnarray}
Given the above conditions, we have the following result.

\begin{lem}\label{lemgeneralvarphibound}
Let $p > 2$, and assume that:
\begin{longlist}[(iii)]
\item[(i)] $\mathbb{E}_{\mathcal{H}_j} [H_j ] = 0$
uniformly for
$j = 1,\ldots,J$,
\item[(ii)] there exists a $u^- > 0$ such that $P (\mathbb
{E}_{\mathcal
{H}_j} [H_j^2 ] \leq u^- ) < 1/7$ uniformly for $j =
1,\ldots,J$,
\item[(iii)] $\mathbb{E} [|H_j|^p ] \leq c_1 < \infty$ uniformly
for $j = 1,\ldots,J$.
\end{longlist}
Then there exist constants $c_{\varphi,1},c_{\varphi,2} > 0$, only
depending on $u^-, c_1$ and $p$, such that\vspace*{-3pt}
\begin{eqnarray*}
&& \mathbb{E} \Biggl[\prod_{j = 1}^{J} \bigl|
\varphi_j^{\mathcal
{H}}(x) \bigr| \Biggr] \lesssim e^{-c_{\varphi,1}x^2 J} +
e^{-\sqrt{J/4}\log
8/7}\qquad \mbox{for $x^2 \leq c_{\varphi,2}$.}
\end{eqnarray*}
\end{lem}

\begin{pf}
Let
\begin{eqnarray*}
&& I(s,x) = \mathbb{E}_{\mathcal{H}_j} \bigl[H_j^2 \bigl(
\bigl(\cos (sxH_j) - \cos(0) \bigr) + \mathrm{i} \bigl(
\sin(sxH_j) \bigr) - \sin(0) \bigr) \bigr].
\end{eqnarray*}
Using a Taylor expansion
%
%
and writing $e^{\mathrm{i}x} = \cos(x) + \mathrm{i}\sin(x)$, we
obtain that
\begin{eqnarray*}
&& \mathbb{E}_{\mathcal{H}_j} \bigl[e^{\mathrm{i}xH_j} \bigr] = 1 -
\mathbb{E}_{\mathcal
{H}_j} \bigl[H_j^2 \bigr]
x^2/2 + x^2\int_0^1(1-s)I(s,x)\,ds.
\end{eqnarray*}
Using the Lipschitz property of $\cos(y)$, $\sin(y)$ and $|e^{\mathrm{i}x}| = 1$, it follows that
%
\begin{eqnarray}
\label{eqlemgeneralvarphibound1}
&& \bigl|I(s,x) \bigr| \leq2\mathbb{E}_{\mathcal{H}_j} \bigl[H_j^2|xh|
+ H_j^2\mathbh{1}\bigl(|H_j|\geq h\bigr) \bigr],\qquad h
> 0.
\end{eqnarray}
For $h > 0$ we have from the Markov inequality
\begin{eqnarray*}
\mathbb{E}_{\mathcal{H}_j} \bigl[H_j^2
\mathbh{1}\bigl(|H_j|\geq h\bigr) \bigr] &\leq& 2 \int_h^{\infty}
x P_{\mathcal{H}_j} \bigl(|H_j| \geq x \bigr)\,dx
+ h^2 P_{\mathcal{H}_j}\bigl(|H_j| \geq h\bigr)
\\
&\leq & 2 h^{-p + 2}\int_0^{\infty}
x^{p-1} P_{\mathcal{H}_j} \bigl(|H_j| \geq x \bigr)\,dx+h^2 P_{\mathcal{H}_j}\bigl(|H_j| \geq h\bigr)
\\
&\leq&  \frac{2+p}{p} h^{-p+2} \mathbb{E}_{\mathcal{H}_j}
\bigl[|H_j|^p \bigr] < 2h^{-p+2}
\mathbb{E}_{\mathcal{H}_j} \bigl[|H_j|^p \bigr].
\end{eqnarray*}
We thus conclude from (\ref{eqlemgeneralvarphibound1}) that
\begin{eqnarray*}
&& \bigl|I(s,x) \bigr|\leq2\mathbb{E}_{\mathcal{H}_j} \bigl[H_j^2
\bigr]|xh| + 4 h^{-p+2} \mathbb{E}_{\mathcal{H}_j} \bigl[|H_j|^p
\bigr].
\end{eqnarray*}
This gives us
%
\begin{eqnarray}
&& \bigl|\mathbb{E}_{\mathcal{H}_j} \bigl[e^{\mathrm{i}x H_j} \bigr] - 1 +
\mathbb{E} _{\mathcal{H}_j} \bigl[H_j^2 \bigr]
x^2/2 \bigr|
\nonumber
\\[-8pt]
\label{eqlemgeneralvarphibound2}
\\[-8pt]
\nonumber
&&\qquad\leq\mathbb {E}_{\mathcal
{H}_j} \bigl[H_j^2
\bigr]h|x|^3 + 2h^{-p+2}x^2 \mathbb{E}_{\mathcal
{H}_j}
\bigl[|H_j|^p \bigr].
\end{eqnarray}
Let $\mathcal{I} =  \{1,\ldots,J \}$, and put $\sigma
_j^{\mathcal{H}} = \mathbb{E} [H_j^2 |\mathcal{H}_j
]$ and $\rho
_j^{\mathcal{H}} = \mathbb{E} [|H_j|^p |\mathcal{H}_j
]$. Consider
\begin{eqnarray*}
&& \rho_{1,J}^{\mathcal{H}} \geq\rho_{2,J}^{\mathcal{H}}
\geq\cdots \geq\rho _{J,J}^{\mathcal{H}},
\end{eqnarray*}
where $\rho_{j,J}^{\mathcal{H}}$ denotes the $j$th largest random variable
for $1 \leq j \leq J$. Let $E_{j}$, $j = 1,\ldots,J$ denote i.i.d.
unit exponential random variables, and denote with $E_{j,J}$ the
$j$th largest. Further, denote with $F_{\rho_j}(\cdot)$ the c.d.f. of
${\rho}_j^{\mathcal{H}}$, $j = 1,\ldots,J$, and with $F_{\rho}(\cdot) =
\min_{1 \leq j \leq J}F_{\rho_j}(\cdot)$. Using the
transformation $-\log (1 - F_{\rho_j}(\cdot) )$, we thus obtain
%
\begin{eqnarray}
&& \bigl\{{\rho}_{j}^{\mathcal{H}} \leq x_j  :
1 \leq j \leq J \bigr\}
\nonumber\\[-8pt]\\[-8pt]\nonumber
&&\qquad \stackrel{d}{=}
\bigl\{E_j \leq -\log\bigl(1 - F_{\rho_j}(x_j)\bigr)  : 1 \leq j \leq J \bigr\}, \qquad x_j \in \mathbb{R},
\end{eqnarray}
which is the well-known Renyi representation; cf. \cite{csoghovbook1993,shorak1986}. In particular, by the
construction of $F_{\rho}(\cdot)$ it follows that
\begin{eqnarray*}
&& P \bigl(\rho_{J/2,J}^{\mathcal{H}} \leq u^+ \bigr) \geq P
(E_{J
/2,J} \leq-\log \bigl(1 - F_{\rho}\bigl(u^+\bigr) \bigr)
\end{eqnarray*}
for $0 \leq u^+ < \infty$. Let $u_{\mathcal{H}}^+ = -\log (1 -
F_{\rho
}(u^+) )$ and $u_{\mathcal{H}}^+(J) = \sqrt{J/2}(u_{\mathcal
{H}}^+ - \log
2)$. We wish to find a $u^+$ such that $u_{\mathcal{H}}^+ > \log2$.
This is
implied by $F_{\rho}(u^+) > 1/2$. We will now construct such an $u^+$. Since
\begin{eqnarray*}
&& c_1 \geq\mathbb{E} \bigl[ \bigl|\rho_j^{\mathcal{H}} \bigr|
\bigr] = \int_0^{\infty} \bigl(1 -
F_{\rho_j}(x) \bigr)\,dx \geq c_2 P \bigl(\rho
_j^{\mathcal{H}} \geq c_2\bigr)\qquad \mbox{for $c_2 > 0$},
\end{eqnarray*}
it follows that $c_1/c_2 \geq1 - P (\rho_j^{\mathcal{H}} < c_2)$. Hence choosing $u^+ = c_2 = 4c_1$, we obtain $F_{\rho_j}(u^+) \geq
3/4$ and hence $F_{\rho}(u^+) \geq3/4$, which leads to $u_{\mathcal{H}}^+
\geq\sqrt{J/2}\log2$. Thus by known properties of exponential
order statistics (cf. \cite{csoghovbook1993,fellervolume2}), we have
\begin{eqnarray*}
P \bigl(E_{J/2,J} \leq u_{\mathcal{H}}^+ \bigr) &= & P \bigl(\sqrt{J
/2}(E_{J/2,J} - \log2) \leq\sqrt{J/2}\bigl(u_{\mathcal{H}}^+ - \log2
\bigr) \bigr)
\\
&=& 1 - P \bigl(\sqrt{J/2}(E_{J/2,J} - \log2) > u_{\mathcal{H}}^+(J)
\bigr)
\\
&\geq & 1 - \mathbb{E} \bigl[e^{\sqrt{J/2}(E_{J
/2,J} - \log2)} \bigr] e^{-u_{\mathcal{H}}^+(J)} \geq1 -
\mathcal {O} \bigl(e^{-\sqrt{J/2}\log2} \bigr),
\end{eqnarray*}
for sufficiently large $J$. We thus conclude that
%
\begin{eqnarray}
\label{eqlemboundchar6}
&& P \bigl(\rho_{J/2,J}^{\mathcal{H}} \leq u^+ \bigr) \geq1
-\mathcal {O} \bigl(e^{-\sqrt{J/2}\log2} \bigr).
\end{eqnarray}
Let us denote this set with $\mathcal{A}^+ =  \{\rho_{J/2,J
}^{\mathcal{H}} \leq u^+  \}$, and put $\mathcal{I}_{\mathcal
{A}}^+ =
 \{1 \leq j \leq J\dvtx \rho_j^{\mathcal{H}} \leq u^+  \}$. Note
that the index set $\mathcal{I}_{\mathcal{A}}^+$ has at least
cardinality $J/2$ given event $\mathcal{A}^+$. For the sake of
simplicity, let us assume that $ |\mathcal{I}_{\mathcal
{A}}^+ | = J/2$, which, as is clear from the arguments below,
has no impact on our results. Let us introduce
\begin{eqnarray*}
&& \sigma_{1,J/2}^{\mathcal{H},\diamond} \geq\sigma_{2,J/2}^{\mathcal{H}
,\diamond}
\geq\cdots \geq\sigma_{J/2,J/2}^{\mathcal{H},\diamond}
\end{eqnarray*}
the order statistics of $\sigma_j^{\mathcal{H}}$ within the index set
$\mathcal{I}_{\mathcal{A}}^+$. This means that $\sigma_{J/2,J
/2}^{\mathcal{H},\diamond}$ is not necessarily the smallest value of
$\sigma
_{j}^{\mathcal{H}}$, $1 \leq j \leq J$. More generally, it holds that
%
\begin{eqnarray}
\label{eqlemboundchar6.5}
&& \sigma_{j,J/2}^{\mathcal{H},\diamond} \geq\sigma_{J/2 +
j,J}^{\mathcal{H}
}, \qquad
j \in\{1,\ldots,J/2\}.
\end{eqnarray}
Now, similar to as before, let $F_{\sigma_j}(\cdot)$ be the c.d.f. of
${\sigma}_j^{\mathcal{H}}$, $j = 1,\ldots,J$, and put $F_{\sigma}(\cdot
) =
\max_{1 \leq j \leq J} F_{\sigma_j}(\cdot)$, $u_{\mathcal{H}}^- =
-\log
 (1 - F_{\sigma}(u^-) )$ and $u_{\mathcal{H}}^-(J) = \sqrt{J
/4}\times\break (\log4/ 3 - u_{\mathcal{H}}^-)$ for some $0 \leq u^- \leq u^+$. We search
for a $u^- >0$ such that $u_{\mathcal{H}}^- < \log7/6$, which is true if
$\max_{1 \leq j \leq J}F_{\sigma_j}(u^-) <1/7$. However, this is
precisely what we demanded in the assumptions. Then proceeding as
before, we have
\begin{eqnarray*}
P \bigl(E_{3J/4,J} \geq u_{\mathcal{H}}^- \bigr) &=& P \bigl(\sqrt{J
/4}(E_{3J/4,J} - \log4/3) \geq\sqrt{J/4}\bigl(u_{\mathcal{H}}^- -
\log4/3\bigr) \bigr)
\\
&=& 1 - P \bigl(\sqrt{J/4}( \log4/3 - E_{3J/4,J}) > u_{\mathcal{H}}^-(J)
\bigr)
\\
&\geq & 1 - \mathbb{E} \bigl[e^{\sqrt
{J/4}(\log
4/3 - E_{3J/4,J})} \bigr] e^{-u_{\mathcal{H}}^-(J)} \\
&\geq & 1 -
\mathcal{O} \bigl(e^{-\sqrt{J/4}\log8/7} \bigr).
\end{eqnarray*}
We thus conclude from (\ref{eqlemboundchar6.5}) and the
construction of $F_{\sigma}(\cdot)$ that
%
\begin{eqnarray}
P \bigl(\sigma_{J/4,J/2}^{\mathcal{H},\diamond} \geq u^-
\bigr) &\geq &  P \bigl(\sigma_{3J/4,J}^{\mathcal{H}} \geq u^- \bigr) \geq P
\bigl(E_{3J/4,J} \geq u_{\mathcal{H}}^- \bigr)
\nonumber
\\[-8pt]
\label{eqlemboundchar7}
\\[-8pt]
\nonumber
&\geq & 1 - \mathcal {O} \bigl(e^{-\sqrt{J/4}\log8/7} \bigr).
\end{eqnarray}
Put $\mathcal{I}_{\mathcal{A}} =  \{j \in\mathcal{I}_{\mathcal
{A}}^+\dvtx  \sigma_{j}^{\mathcal{H},\diamond} \geq u^-  \}$. Combining
(\ref{eqlemboundchar6}) and (\ref{eqlemboundchar7}) we obtain
%
\begin{eqnarray}
\label{eqlemboundchar8}
&& P \bigl( \bigl\{\sigma_{J/4,J/2}^{\mathcal{H},\diamond} \geq u^-
\bigr\} \cap \bigl\{\rho_{J/2,J}^{\mathcal{H}} \leq u^+ \bigr\} \bigr)
\geq1 - \mathcal{O} \bigl(e^{-\sqrt{J/4}\log8/7} \bigr).
\end{eqnarray}
We denote this set with $\mathcal{A}=  \{\sigma
_{J/4,J/2}^{\mathcal{H}
,\diamond} \geq u^-  \} \cap \{\rho_{J/2,J}^{\mathcal{H}}
\leq u^+ \}$. Also note that by the (conditional) Lyapunov
inequality, we have
%
\begin{eqnarray}
&& \rho_{j,J}^{\mathcal{H}} \geq \bigl(\sigma_{j,J}^{\mathcal
{H}}
\bigr)^{p/2}.
\end{eqnarray}
Note that $|\mathcal{I}_{\mathcal{A}}| \geq J/4$ on the event
$\mathcal{A}
$, and, by the above, we get
%
\begin{eqnarray}
&& \rho_{j}^{\mathcal{H}} \leq u^+ \quad\mbox{and}\quad u^- \leq\sigma
_{j}^{\mathcal{H}} \leq \bigl(u^+ \bigr)^{p/2}\qquad \mbox{for
$j \in \mathcal{I}_{\mathcal{A}}$.}
\end{eqnarray}
Using (\ref{eqlemgeneralvarphibound2}), this implies that for
every $j \in\mathcal{I}_{\mathcal{A}}$, we have
\begin{eqnarray*}
&&\bigl|\mathbb{E}_{|\mathcal{H}_j} \bigl[e^{\mathrm{i}x H_j} \bigr] - 1 + \mathbb{E}
_{|\mathcal{H}_j} \bigl[H_j^2 \bigr] x^2/2 \bigr|
\leq\bigl(u^+\bigr)^{p/2} h|x|^3 + 2h^{-p+2} u^+
x^2.
\end{eqnarray*}
%
Hence, if $(u^+)^{p/2} x^2/2 < 1$ and $h = x^{-1/(p-1)}$, we conclude
from the above and the triangle inequality that for $j \in\mathcal
{I}_{\mathcal{A}}$,
%
\begin{eqnarray}
&& \bigl|\varphi_j^{\mathcal{H}}(x) \bigr| < 1 - u^-
x^2/2 + \bigl(2u^+ + \bigl(u^+\bigr)^{p/2} \bigr)
|x|^{2 + \delta(p)}
\nonumber
\\[-8pt]
\label{eqlemboundchar9}
\\[-8pt]
\eqntext{\mbox{for $\bigl(u^+\bigr)^{p/2} x^2/2 <
1$,}}
\end{eqnarray}
where $\delta(p) = (p-2)/(p-1) > 0$. Since $0 < u^-, u^+ < \infty$ and $\delta(p) > 0$, there exist
absolute constants $0 < c_{\varphi,1}, c_{\varphi,2}$, chosen
sufficiently small, such that
%
\begin{eqnarray}
\label{eqlemboundchar9.2}
&&\qquad u(x) \stackrel{\mathrm{def}} {=} u^- x^2 /2 - \bigl(2u^+ +
\bigl(u^+\bigr)^{p/2} \bigr) |x|^{2 + \delta(p)} \geq8 c_{\varphi,1}
x^2\qquad \mbox{for $x^2 \leq c_{\varphi,2}$}.\hspace*{-12pt}
\end{eqnarray}
Next, observe that\vspace*{-3pt}
%
\begin{eqnarray}
\nonumber
\mathbb{E} \Biggl[ \Biggl|\prod_{j = 1}^{J}
\varphi_j^{\mathcal
{H}}(x) \Biggr| \Biggr] & =& \mathbb{E} \Biggl[ \Biggl|\prod
_{j = 1}^{J}\varphi _j^{\mathcal{H}
}(x)
\Biggr| \bigl(\mathbh{1}(\mathcal{A}) + \mathbh{1}\bigl(\mathcal {A}^c
\bigr) \bigr) \Biggr] \\[-3pt]
\label{eqlemboundchar9.1}
&\leq &  P \bigl(\mathcal{A}^c \bigr) + \mathbb{E}
\biggl[ \biggl|\prod_{j \in\mathcal{I}}\varphi_j^{\mathcal{H}}(x)
\biggr|\mathbh {1}(\mathcal {A}) \biggr]
\\[-3pt]
\nonumber
& \leq &  P \bigl(\mathcal{A}^c \bigr) + \mathbb{E} \biggl[ \biggl|\prod
_{j \in\mathcal{I}_{\mathcal{A}}}\varphi_j^{\mathcal
{H}}(x) \biggr|
\mathbh{1}(\mathcal{A}) \biggr].
\end{eqnarray}
Moreover, using (\ref{eqlemboundchar9}) and (\ref{eqlemboundchar9.2}) and since $|\mathcal
{I}_{\mathcal{A}}| = J/4$ on $\mathcal{A}$, it follows that
\begin{eqnarray*}
\mathbb{E} \biggl[ \biggl|\prod_{j \in\mathcal{I}_{\mathcal
{A}}}\varphi
_j^{\mathcal{H}}(x) \biggr|\mathbh{1}(\mathcal{A}) \biggr] &\leq &
\mathbb{E} \biggl[ \biggl|\prod_{j \in\mathcal{I}_{\mathcal{A}}} \bigl(1 - u(x) \bigr)
\biggr|\mathbh{1}(\mathcal{A}) \biggr]
\\
&\leq & \mathbb{E} \biggl[\prod_{j
\in\mathcal
{I}_{\mathcal{A}}}e^{- u(x)}
\mathbh{1}(\mathcal{A}) \biggr]\leq e^{-u(x)
J/8} \leq e^{-c_{\varphi,1}Jx^2}.
\end{eqnarray*}
Hence we conclude from the above and (\ref{eqlemboundchar8}) that
\begin{eqnarray}
\nonumber
&&\mathbb{E} \biggl[ \biggl|\prod_{j \in\mathcal{I}}\varphi
_j^{\mathcal{H}
}(x) \biggr| \biggr] \lesssim e^{-c_{\varphi,1}x^2 J} +
e^{-\sqrt
{J/4}\log8/7}\qquad \mbox{for $x^2 \leq c_{\varphi,2}$,}
\end{eqnarray}
which yields the claim.
\end{pf}

\begin{lem}\label{lemsigexpressionsrelations}
Grant Assumption~\ref{assdependence}. Then $\sum_{k = 1}^{\infty}k
 |\mathbb{E} [X_0 X_k ] | < \infty$ and $\widehat
{\sigma
}_m^2 = s_m^2/2 + \mathcal{O} (m^{-1} )$. Moreover, we have
$\widehat{\sigma}_l^2 = s_m^2/2 + \mbox{\scriptsize$\mathcal
{O}$} (1 )$ as $l \to m$.
\end{lem}

\begin{pf}
Since $\mathbb{E} [X_k |\mathcal{E}_0 ] = \mathbb
{E} [X_k -
X_k^* |\mathcal{E}_0 ]$, the Cauchy--Schwarz and Jensen
inequalities imply
\begin{eqnarray*}
\sum_{k = 0}^{\infty} \bigl|\mathbb{E}
[X_0X_k ] \bigr| &\leq &  \|X_0 \|_2
\sum_{k = 0}^{\infty} \bigl\|\mathbb{E}
[X_k |\mathcal{E}_0 ] \bigr\|_2 \leq
\|X_0 \|_2\sum_{k =
0}^{\infty}
\bigl\|X_k-X_k^* \bigr\|_2
\\
&\leq&  \|X_0 \| _2\sum_{k = 1}^{\infty}k^2
\bigl\|X_k-X_k' \bigr\|_2 < \infty.
\end{eqnarray*}
The decomposition $\widehat{\sigma}_m^2 = s_m^2/2 + \mathcal{O}
(m^{-1} )$ now follows from (\ref{decompsigma}). Claim $\widehat
{\sigma}_l^2 = s_m^2/2 + \mathcal{O}(1)$ as $l \to m$ readily follows from the previous computations.
\end{pf}

\begin{lem}\label{lemsigexpansion}
Grant Assumption~\ref{assdependence}. Then:
\begin{longlist}[(iii)]
\item[(i)] $\|\sigma_{j|m}^2 - \sigma_j^2 \|_{p/2}
\lesssim \|\sigma_{j|m}^2 - \widehat{\sigma}_m^2 \|
_{p/2}+m^{-1} \lesssim m^{-1}$ for $1 \leq j \leq N$,
\item[(ii)] $\sigma_j^2 = \widehat{\sigma}_m^2 + \mathcal{O}
(m^{-1} )$ for $1 \leq j \leq N$,
\item[(iii)] $ \|\sigma_{|m}^2 - \overline{\sigma}_m^2 \|
_{p/2} \lesssim n^{-1} N^{2/p}$.
\end{longlist}
\end{lem}

\begin{pf}
We first show (i). Without\vspace*{1pt} loss of generality, we may assume $j = 1$,
since $m \thicksim m'$. To lighten the notation, we use $R_1 =
R_1^{(1)}$. We will first establish that $ \|\sigma_{j|m}^2 -
\widehat{\sigma}_m^2 \|_{p/2} \lesssim m^{-1}$. We have that
\begin{eqnarray*}
&& 2m \bigl(\sigma_{1|m}^2 - \widehat{\sigma}_m^2
\bigr)\\
&&\qquad= \mathbb {E}_{\mathbb{F}
_m} \Biggl[ \Biggl(\sum
_{k = 1}^m \bigl(X_k^{(**)} +
\bigl(X_k - X_k^{(**)}\bigr) -
\mathbb{E}_{\mathbb{F}_m}[X_k] \bigr) + R_1
\Biggr)^2 \Biggr] - 2m \widehat{\sigma}_m^2.
\end{eqnarray*}
By squaring out the first expression, we obtain a sum of square terms
and a sum of mixed terms. Let us first treat the mixed terms, which are
\begin{eqnarray*}
&& 2\sum_{k = 1}^m \sum
_{l = 1}^m \mathbb{E}_{\mathbb{F}_m} \bigl[X_k^{(**)}\bigl(X_l - X_l^{(**)}\bigr)
+ X_k^{(**)}\E_{\FF_m}[X_l] + \E_{\FF_m}[X_k]\bigl(X_l - X_l^{(**)}\bigr)\bigr]
\\[-2pt]
&&\quad{}+ 2 \sum_{k
= 1}^m \mathbb{E}_{\mathbb{F}_m}
\bigl[R_1 X_k^{(**)} +R_1
\bigl(X_k - X_k^{(**)}\bigr) + R_1
\mathbb{E}_{\mathbb{F}_m}[X_k] \bigr]
\\[-2pt]
&&\qquad=  I_m + \mathit{II}_m +\mathit{III}_m + \mathit{IV}_m +
V_m + \mathit{VI}_m.
\end{eqnarray*}
We will handle all these terms separately.
\begin{longlist}[{}]
\item[\textit{Case} $I_m$:] We have
\begin{eqnarray*}
I_m/2 &=& \sum_{l = 1}^m \sum
_{k = l}^m(\cdots) + \sum
_{l = 1}^m \sum_{k = 1}^{l-1}(\cdots)
\\[-2pt]
&=& \sum_{l = 1}^m \sum
_{k = l}^m \mathbb {E}_{\mathbb{F}
_m} \bigl[
\bigl(X_l - X_l^{(**)}\bigr)\mathbb{E}
\bigl[X_k^{(**)} |\sigma \bigl(\mathbb{F}_m,
\mathcal{E}_l,\mathcal{E}_l^{(**)} \bigr)
\bigr] \bigr]
\\[-2pt]
&&{}+ \sum_{l = 1}^m \sum
_{k = 1}^{l-1}\mathbb{E}_{\mathbb{F}_m}
\bigl[X_k^{(**)}\bigl(X_l -
X_l^{(**)}\bigr) \bigr].
\end{eqnarray*}
Since\vspace*{-3pt}
\begin{eqnarray*}
&& \mathbb{E} \bigl[X_k^{(**)} |\sigma \bigl(
\mathbb{F}_m,\mathcal {E}_l,\mathcal
{E}_l^{(**)} \bigr) \bigr] \stackrel{d} {=} \mathbb{E}
[X_k |\mathcal{E}_l ] = \mathbb{E}
\bigl[X_k - X_k^{(k-l,*)} |\mathcal
{E}_l \bigr],
\end{eqnarray*}
%
the Cauchy--Schwarz (with respect to $\mathbb{E}_{\mathbb{F}_m}$) and Jensen
inequalities thus yield
\begin{eqnarray*}
\|I_m \|_{p/2} &\leq & 2\sum_{l = 1}^m
\sum_{k = l}^m \bigl\|X_l -
X_l^{(**)} \bigr\|_{p} \bigl\|X_k -
X_k^{(k-l,*)} \bigr\|_p \\
&&{}+ 2\sum
_{l = 1}^m \sum_{k = 1}^{l-1}
\bigl\|X_k^{(**)} \bigr\|_{p} \bigl\|X_l -
X_l^{(**)} \bigr\|_{p}
\\
& \leq & 2 \Biggl(\sum_{l = 1}^{\infty
}l
\bigl\|X_l-X_l' \bigr\|_p
\Biggr)^2 + 2\sum_{l = 1}^{\infty
}l^2
\bigl\|X_l-X_l' \bigr\|_p
\|X_1 \|_p < \infty.
\end{eqnarray*}
\item[\textit{Case} $\mathit{II}_m$:] Since $\mathbb{E}_{\mathbb{F}_m}[X_k^{(**)}] =\mathbb
{E}[X_k] = 0$ ($k \leq m$) it follows
that $\mathit{II}_m = 0$.
\item[\textit{Case} $\mathit{III}_m$:] It follows via the Jensen and triangle inequalities that\vspace*{-2pt}
%
\begin{eqnarray}
 \bigl\|\mathbb{E} [X_l |\mathbb{F}_m ]
\bigr\|_p &= & \bigl\|\mathbb{E} \bigl[X_l - X_l^{(**)}
|\mathbb{F}_m \bigr] \bigr\|_p
\nonumber
\\[-10pt]
\label{eqlemsigexpansion3}\\[-10pt]
\nonumber
 &\leq &  \bigl\|X_l -
X_l^{(**)}\bigr\|_p \leq\sum
_{j = l}^{\infty} \bigl\|X_j -
X_j' \bigr\|_p.
\end{eqnarray}
The Cauchy--Schwarz (with respect to $\mathbb{E}_{\mathbb{F}_m}$) and Jensen
inequalities then give\vspace*{-3pt}
\begin{eqnarray*}
&& \|\mathit{III}_m \|_{p/2} \leq2 \Biggl(\sum
_{l = 1}^m \sum_{j =
l}^{\infty}
\bigl\|X_j - X_j' \bigr\|_p
\Biggr)^2 \leq2 \Biggl(\sum_{l = 1}^{\infty}
l \bigl\|X_l - X_l' \bigr\|_p
\Biggr)^2 < \infty.
\end{eqnarray*}
\item[\textit{Case} $\mathit{IV}_m$:] Note that $X_k^{(**)}$ and $X_l^{(l-k,*)}$ are independent
for $1 \leq k \leq m$ and $m+1\leq l \leq2m$. Hence since $\mathbb
{E}_{\mathbb{F}
_m} [X_k^{(**)}  ] = 0$, we have\vspace*{-3pt}
\begin{eqnarray*}
\sum_{k = 1}^m \mathbb{E}_{\mathbb{F}_m}
\bigl[R_1 X_k^{(**)} \bigr] &=& \sum
_{k =
1}^m\sum_{l = m+1}^{2m}
\mathbb{E}_{\mathbb{F}_m} \bigl[X_k^{(**)}X_l
\bigr]
\\[-3pt]
&=& \sum_{k = 1}^m\sum
_{l = m+1}^{2m} \mathbb{E}_{\mathbb{F}_m}
\bigl[X_k^{(**)}\bigl(X_l -
X_l^{(l-k,*)} + X_l^{(l-k,*)}\bigr)
\bigr]
\\[-3pt]
&=& \sum_{k = 1}^m\sum
_{l =
m+1}^{2m} \mathbb{E}_{\mathbb{F}_m}
\bigl[X_k^{(**)}\bigl(X_l -
X_l^{(l-k,*)}\bigr) \bigr].
\end{eqnarray*}
The Cauchy--Schwarz (with respect to $\mathbb{E}_{\mathbb{F}_m}$) and Jensen
inequalitities then yield\vspace*{-3pt}
\begin{eqnarray*}
\|\mathit{IV}_m \|_{p/2} &\leq & 2 \sum_{k = 1}^m
\sum_{l = m+1}^{2m} \bigl\|X_k^{(**)}
\bigr\|_{p} \bigl\|X_l - X_l^{(l-k,*)}
\bigr\|_{p}
\\[-3pt]
&\leq & 2\sum_{k = 1}^m\sum
_{l = m+1}^{2m} \sum_{j = l-k}^{\infty
}
\bigl\|X_j - X_j' \bigr\|_{p}
\|X_1 \|_{p}
\\[-3pt]
& \leq &2\sum_{k = 1}^m\sum
_{j = m-k}^{\infty} (j-m+k) \bigl\|X_j -
X_j' \bigr\| _{p} \|X_1
\|_{p}
\\[-3pt]
&\leq & 2 \sum_{j = 1}^{\infty} j^2
\bigl\|X_j - X_j' \bigr\|_{p}
\|X_1 \|_{p} < \infty.
\end{eqnarray*}
\item[\textit{Case} $V_m$:] The Cauchy--Schwarz (with respect to $\mathbb{E}_{\mathbb
{F}_m}$) and
Jensen inequalities yield\vspace*{-9pt}
\begin{eqnarray*}
 \|V_m \|_{p/2} &\leq & 2 \sum_{k = 1}^m
\bigl\|X_k - X_k^{(**)} \bigr\|_p
\|R_1 \|_p \\[-3pt]
&\leq & \sum_{l = 1}^{\infty
}l^2
\bigl\|X_l - X_l' \bigr\|_p
\|R_1 \|_p < \infty,
\end{eqnarray*}
since\vadjust{\goodbreak} $ \|R_1 \|_p < \infty$ by
Lemma~\ref{lemboundR1}.\vadjust{\goodbreak}

\item[\textit{Case} $\mathit{VI}_m$:] Proceeding as above and using (\ref{eqlemsigexpansion3}), we get
$ \|\mathit{VI}_m \|_{p/2} < \infty$. It thus
remains to
deal with the squared terms, which are
\begin{eqnarray*}
&&\sum_{k = 1}^m\sum
_{l = 1}^m\mathbb{E}_{\mathbb{F}_m}
\bigl[X_k^{(**)}X_l^{(**)} +
\bigl(X_k - X_k^{(**)}\bigr)
\bigl(X_l - X_l^{(**)}\bigr)\\
&&\hspace*{119pt}\qquad {}+
\mathbb{E}_{\mathbb
{F}_m}[X_k]\mathbb{E}_{\mathbb{F}
_m}[X_l]
\bigr] + \mathbb{E}_{\mathbb{F}_m} \bigl[R_1^2 \bigr]
\\
&&\qquad= 2m\widehat {\sigma}_m^2 + \mathit{VII}_m +
\mathit{VIII}_m + \mathit{IX}_m.
\end{eqnarray*}
However, using the results from the previous computations and
Lemma~\ref{lemboundR1}, one readily deduces that
%
\begin{eqnarray}
&& \quad\|\mathit{VII}_m \|_{p/2} < \infty,\qquad \| \mathit{VIII}_m
\|_{p/2} < \infty,\qquad \|\mathit{IX}_m \| _{p/2} < \infty.
\end{eqnarray}
Piecing everything together, we have established that $ \|\sigma
_{j|m}^2 - \widehat{\sigma}_m^2 \|_{p/2} \lesssim m^{-1}$.
However, from the above arguments one readily deduces that $\sigma_j^2
= \widehat{\sigma}_m^2 + \mathcal{O} (m^{-1} )$, and hence
(i) and
(ii) follow. We now treat (iii). Since $ \{Y_j^{(1)} \}_{1
\leq j \leq N}$ is an independent sequence under $P_{\mathbb{F}_m}$,
we have
%
\begin{eqnarray}
\label{eqlemsigexpansion6}
&& \sigma_{|m}^2 = N^{-1} \sum
_{j = 1}^N \sigma_{j|m}^2.
\end{eqnarray}
Let $\mathcal{I} = \{1,3,5,\ldots\}$ and $\mathcal{J} = \{2,4,6,\ldots\}$
such that $\mathcal{I} \cup\mathcal{J} = \{1,2,\ldots,N\}$. Then
\begin{eqnarray*}
&&\Biggl\|\sum_{j = 1}^N \sigma_{j|m}^2
- \sigma_j^2 \Biggr\|_{p/2} \leq \biggl\|\sum
_{j \in\mathcal{I}} \sigma_{j|m}^2 - \sigma
_j^2 \biggr\|_{p/2} + \biggl\|\sum
_{j \in\mathcal{J}} \sigma _{j|m}^2 -
\sigma_j^2 \biggr\|_{p/2}.
\end{eqnarray*}
Note that $\{\sigma_{j|m}^2\}_{j \in\mathcal{I}}$ is a sequence of
independent random variables, and the same is true for $\{\sigma
_{j|m}^2\}_{j \in\mathcal{J}}$. Then by Lemma~\ref{lemwuoriginal},
it follows that
\begin{eqnarray*}
&& \Biggl\|\sum_{j = 1}^N \sigma_{j|m}^2
- \sigma_j^2 \Biggr\|_{p/2} \lesssim
N^{2/p} \bigl\|\sigma_{j|m}^2 - \sigma_j^2
\bigr\|_{p/2}\qquad \mbox{for $p \in(2,3]$,}
\end{eqnarray*}
which by (i) is of the magnitude $\mathcal{O} (N^{2/p}
m^{-1} )$.
Hence we conclude from (\ref{eqlemsigexpansion6}) that
\begin{eqnarray*}
&& \bigl\|\sigma_{|m}^2 - \overline{\sigma}_m^2
\bigr\|_{p/2} \lesssim n^{-1} N^{2/p}.
\end{eqnarray*}
\end{longlist}
\upqed\end{pf}

\begin{cor}\label{corsigreplaceissmall}
Grant Assumption~\ref{assdependence}. Let $\mathcal{B} =  \{
\sigma_{|m}^2 \geq s_m^2/4  \}$. Then
\begin{eqnarray*}
&& P \bigl(\mathcal{B}^c \bigr)\lesssim n^{-p/2} N.
\end{eqnarray*}
\end{cor}

\begin{pf}
By Markov's inequality, Lemma~\ref{lemsigexpressionsrelations} and
Lemma~\ref{lemsigexpansion}, it follows that for large enough $m$
\begin{eqnarray*}
&& P \bigl(\mathcal{B}^c \bigr) \leq P \bigl( \bigl|\sigma_{|m}^2
- \overline{\sigma}_m^2 \bigr| \geq s_m^2/4
- \mathcal {O}\bigl(m^{-1}\bigr) \bigr) \lesssim s_m^{-p/2}n^{-p/2}
N \lesssim n^{-p/2} N
\end{eqnarray*}
since $s_m^2 > 0$ by Assumption~\ref{assdependence}(iii).
\end{pf}

\subsubsection{Ideal metrics and applications}\label{seczolo}
The aim of this section is to give a proof for the inequality
%
\begin{eqnarray}
\label{eqphiboundcor}
&& \bigl\|\varphi_i(t) -e^{-\sigma_{i|m}^2 t^2/2} \bigr\|_1
\lesssim |t|^p m^{-p/2 + 1},\qquad i \in\{1,N\}
\end{eqnarray}
in Corollary~\ref{corzolo}. We will achieve this by employing ideal
metrics. Let $s > 0$. Then we can represent $s$ as $s = \mathfrak{m}+
\alpha$,
where $[s] = \mathfrak{m}$ denotes the integer part, and $0 < \alpha
\leq1$.
Let $\mathfrak{F}_s$ be the class of all real-valued functions $f$,
such that the $\mathfrak{m}$th derivative exists and satisfies
%
\begin{eqnarray}
&& \bigl|f^{(\mathfrak{m})}(x) - f^{(\mathfrak{m})}(y) \bigr| \leq |x - y |^{\alpha}.
\end{eqnarray}
Note that since $\cos(y), \sin(y)$ are bounded in absolute value and
are Lipschitz continuous, it follows that up to some finite constant
$c(\alpha)>0$ we have $\sin(y),\cos(y) \in\mathfrak{F}_s$ for any
$s > 0$. 
As already mentioned in step~2 of the proof of part ${\mathbf A}$, we will
make use of some special ideal-metrics $\zeta_s$ (Zolotarev metric).
For two probability measures $P,Q$, the metric $\zeta_s$ is defined as
\begin{eqnarray*}
&& \zeta_s(P,Q) = \sup \biggl\{ \biggl|\int f(x) (P - Q) (dx) \biggr|\dvtx  f \in
\mathfrak{F}_s \biggr\}.
\end{eqnarray*}
The metric $\zeta_s(P,Q)$ has the nice property of homogeneity. For
random variables $X,Y$, induced probability measures $P_{cX}, P_{cY}$
and constant $c > 0$, this means that $\zeta_s(P_{cX},P_{cY}) = |c|^s
\zeta_s(P_X,P_Y)$. We require some further notation. For $1 \leq j
\leq N$, put
\begin{eqnarray*}
&& S_{j|m} = \frac{1}{\sqrt{2m}}\sum_{k = (2j-2)m + 1}^{(2j-1)m}
X_k\quad \mbox{and}\quad S_{j|m}^{(**)} = \frac{1}{\sqrt{2 m}}
\sum_{k
= (2j-2)m + 1}^{(2j-1)m} X_k^{((2j - 2)m, **)}.
\end{eqnarray*}
Note that $S_{j|m}^{(**)}$ is independent of $\mathbb{F}_m$, and hence
$\widehat{\sigma}_m^2 = \mathbb{E}_{\mathbb{F}_m} [
(S_{j|m}^{(**)}
)^2 ]$. Let $ \{Z_{j} \}_{1 \leq j \leq N}$ be a
sequence of zero mean, standard i.i.d. Gaussian random variables. In
addition, let
\begin{eqnarray*}
&& \eta_{j|m}^2 = \frac{1}{\sqrt{2m}} \mathbb{E}_{\mathbb{F}_m}
\bigl[\bigl(Y_j^{(1)}\bigr)^2 \bigr]/\widehat{
\sigma}_m^2 = \sigma_{j|m}^2/
\widehat {\sigma}_m^2 \qquad \mbox{for $1 \leq j \leq N-1$},
\end{eqnarray*}
and $\eta_{m'|m}^2 = \sigma_{m'|m}^2/\widehat{\sigma}_{m'}^2$ for
$j = N$. As\vspace*{1pt} first step toward (\ref{eqphiboundcor}), we have the
following result.

\begin{lem}\label{lemtrunccharcosandsin}
Grant Assumption~\ref{assdependence}. Then for $f(x) \in \{\cos
(x), \sin(x) \}$, it~holds that
\begin{eqnarray*}
&&\bigl\|\mathbb{E}_{\mathbb{F}_m} \bigl[f\bigl(x(2m)^{-1/2}Y_j^{(1)}
\bigr) - f\bigl(xS_{j|m}^{(**)}\eta_{j|m}\bigr)
\bigr] \bigr\|_1 \lesssim m^{-p/2+1}|x|^p \qquad\mbox{if $p
\in(2,3]$},
\end{eqnarray*}
where $j = 1,\ldots,N$.
\end{lem}

\begin{pf}
To lighten the notation, we use $Y_j = (2m)^{-1/2} Y_j^{(1)}$ and $S_j
= S_{j|m}^{(**)}$ in the following. Using Taylor expansion, we have
%
\begin{eqnarray}
\label{eqlemtrunccharcosandsin2}
&&
\qquad f(y) = f(0) + yf'(0) + y^2/2f''(0)
+ y^2\int_0^1(1-t)
\bigl(f''(ty)-f''(0)
\bigr)\,dt.\hspace*{-10pt}
\end{eqnarray}
Note $\mathbb{E}_{\mathbb{F}_m} [Y_j ] = 0$ and $\mathbb
{E}_{\mathbb{F}_m}
[S_{j}\eta_{j|m} ] = \eta_{j|m}\mathbb{E} [S_{j} ] = 0$.
Moreover, since $\sigma_{j|m}^2 = \eta_{j|m}^2\mathbb{E}_{\mathbb
{F}_m}
[S_j^2 ]$ by construction, we obtain from (\ref
{eqlemtrunccharcosandsin2}) that
\begin{eqnarray*}
&&\mathbb{E}_{\mathbb{F}_m} \bigl[f(xY_j) - f(xS_{j}
\eta_{j|m}) \bigr]
\\
&&\qquad= x^2\int_0^1(1-t)
\mathbb{E}_{\mathbb{F}_m} \bigl[Y_j^2
\bigl(f''(txY_j)-f''(0)
\bigr)\\
&&\qquad\qquad{}\hspace*{68pt}- (S_{j}\eta_{j|m})^2
\bigl(f''(txS_{j}\eta_{j|m})-f''(0)
\bigr) \bigr] \,dt
\\
&&\qquad\stackrel{\mathrm{def}}{=} x^2I_m(x).
\end{eqnarray*}
We have
\begin{eqnarray*}
&& Y_j^2 - S_{j}^2
\eta_{j|m}^2 = Y_j^2 -
S_{j}^2 + S_{j}^2\widehat {
\sigma}_m^{-2} \bigl(\sigma_{j|m}^2 -
\widehat{\sigma}_m^{2} \bigr),
\end{eqnarray*}
where we recall that $S_j$ and $\sigma_{j|m}$ are independent. Using
the Jensen, triangle and H\"{o}lder inequalities and $f/2 \in\mathfrak
{F}_{p}$, it follows that
\begin{eqnarray*}
&& \bigl\|\mathbb{E}_{\mathbb{F}_m} \bigl[ \bigl(Y_j^2 -
(S_{j}\eta _{j|m})^2 \bigr)
\bigl(f''(txY_j) -f''(0)
\bigr) \bigr] \bigr\|_1 \\
 &&\qquad\leq  2 \bigl\| \bigl(Y_j^2 -
(S_{j}\eta_{j|m})^2 \bigr) |txY_j
|^{p-2} \bigr\|_1
\\
&&\qquad\leq  2 \bigl\|Y_j^2 - (S_{j}
\eta_{j|m})^2 \bigr\| _{p/2} \bigl\||txY_j|^{p-2}
\bigr\|_{p/(p-2)}\\
&&\qquad \lesssim2 \bigl\| |txY_j|^{p-2}
\bigr\|_{p/(p-2)}
\\
&&\qquad\quad{}\times\bigl( \|Y_j-S_{j} \|_p
\|Y_j + S_{j} \|_{p} + \bigl\|S_j^2
\bigr\|_{p/2} \bigl\|\sigma_{j|m}^2 - \widehat{\sigma
}_m^{2} \bigr\|_{p/2} \bigr),
\end{eqnarray*}
where we used that $\widehat{\sigma}_m > 0$ for large enough $m$. By
Lemmas~\ref{lemboundR1} and~\ref{lemsigexpansion}, this is
of magnitude $\mathcal{O} (m^{-1/2}|tx|^{p-2} )$. Hence by adding
and subtracting $f''(txY_j)$ and using similar arguments as before, we
obtain from the above
%
\begin{eqnarray}
\nonumber
&& x^2 \bigl\|I_m(x) \bigr\|_{1}
\\
&&\qquad\lesssim   m^{-1/2}|x|^p + x^2\int
_0^1(1-t) \bigl\|\mathbb{E}_{\mathbb{F}_m}
\bigl[S_{j}^2\eta _{j|m}^2
\bigl(f''(txY_j)-f''(txS_{j})
\bigr) \bigr] \bigr\|_1 \, dt
\nonumber
\\[-8pt]
\label{eqlemtrunccharcosandsin3}
\\[-8pt]
\nonumber
&&\qquad
\lesssim   m^{-1/2}|x|^p + |x|^p
\bigl\|S_{j}^2 \bigr\|_{p/2} \bigl\|\eta _{j|m}^2
\bigr\|_{p/2} \bigl\||S_{j} - Y_j|^{p-2}
\bigr\|_{p/(p-2)}
\\
\nonumber
&&\qquad\lesssim   m^{-1/2}|x|^p + m^{-(p-2)/2}|x|^p
\lesssim m^{-p/2 + 1}|x|^p,
\end{eqnarray}
where we use that $S_j$ and $\eta_{j|m} = \widehat{\sigma
}_{j|m}/\widehat{\sigma}_m$ are independent. This gives the desired result.
\end{pf}

As next step toward (\ref{eqphiboundcor}), we have the following.

\begin{lem}\label{lemapproxchar}
Grant Assumption~\ref{assdependence}. Then for $f(x) \in \{\cos
(x), \sin(x) \}$, it holds that
\begin{eqnarray*}
&& \bigl\|\mathbb{E}_{\mathbb{F}_m} \bigl[f(xZ_j\sigma_{j|m}) -
f\bigl(xS_{j|m}^{(**)}\eta_{j|m}\bigr) \bigr]
\bigr\|_1 \lesssim|x|^p m^{-p/2
+ 1}\qquad \mbox{if $p
\in(2,3]$},
\end{eqnarray*}
where $j = 1,\ldots,N$.
\end{lem}

\begin{pf}
To increase the readability, we use the abbreviations $\widehat{\sigma
} = \widehat{\sigma}_m$ and $S_j = S_{j|m}^{(**)}$ in the following.\vspace*{1pt}
The main objective is to transfer the problem to the setup in
\cite{dedeckerminimaldistance2009} and apply the corresponding results.
To this end, we first perform some necessary preparatory computations.
We have that
%
\begin{eqnarray}
&& \sum_{k>l}\bigl\|\mathbb{E} [X_k |
\mathcal{E}_0 ] \bigr\|_p \leq\sum
_{k>l} \bigl\|X_k - X_k^* \bigr\|_p
\leq\sum_{k>l}k \bigl\| X_k -
X_k' \bigr\|_p \to 0
\end{eqnarray}
as $l \to\infty$, hence it follows that
%
\begin{eqnarray}
\label{eqlemapproxchar1}
&& \sum_{k = 0}^{m}\mathbb{E}
[X_k |\mathcal{E}_0 ]\qquad \mbox{converges in $\|\cdot
\|_p$}.
\end{eqnarray}
Next, note that
%
\begin{eqnarray}
&& 2m \widehat{\sigma}^2 = \mathbb{E}_{\mathcal{E}_0} \Biggl[ \Biggl(
\sum_{k
= 1}^m X_k^*
\Biggr)^2 \Biggr].
\end{eqnarray}
Using exactly the same arguments as in the proof of
Lemma~\ref{lemsigexpansion} (the present situation is much
simpler), we get that
\begin{eqnarray*}
&& \Biggl\|\mathbb{E}_{\mathcal{E}_0} \Biggl[ \Biggl(\sum_{k = 1}^m
X_k \Biggr)^2\Biggr] - 2m\widehat{\sigma}^2
\Biggr\|_{p/2} < \infty.
\end{eqnarray*}
We thus obtain that
%
\begin{eqnarray}
\label{eqlemapproxchar2}
&& \sum_{m = 1}^{\infty}
m^{-3/2} \Biggl\|\mathbb{E}_{\mathcal
{E}_0} \Biggl[ \Biggl(\sum
_{k = 1}^m X_k \Biggr)^2
\Biggr] - 2m\widehat{\sigma}^2 \Biggr\|_{p/2} \lesssim\sum
_{m = 1}^{\infty} m^{-3/2} < \infty.
\end{eqnarray}
We will now treat the cases $p \in(2,3)$ and $p = 3$ separately.
\begin{longlist}[{}]
\item[\textit{Case} $p \in(2,3)$:] Since $\cos(y)/2, \sin(y)/2 \in
\mathfrak
{F}_p$, the homogeneity of order $p$ implies that
\begin{eqnarray*}
\bigl\|\mathbb{E}_{\mathbb{F}_m} \bigl[f (xS_{j}\eta _{j|m}
) - f (xZ_{j}\widehat{\sigma}\eta_{j|m} ) \bigr]
\bigr\|_1  &\leq &  2 \bigl\|\zeta_{p} (P_{xS_{j}\eta_{j|m}|\mathbb{F}_m},
P_{xZ_{j}\widehat{\sigma}\eta_{j|m}|\mathbb{F}_m} ) \bigr\| _1
\\
&=& 2 \bigl\||x|^p|\eta_{j|m}|^p\zeta_{p}
(P_{S_{j}|\mathbb{F}_m}, P_{Z_{j}\widehat{\sigma}|\mathbb{F}_m} ) \bigr\|_1.
\end{eqnarray*}
Since $S_j,Z_j$ are independent of $\mathbb{F}_m$, we have $\zeta
_{p}
(P_{S_{j}|\mathbb{F}_m}, P_{Z_{j}\widehat{\sigma}|\mathbb
{F}_m} ) = \zeta
_{p} (P_{S_{j}}, P_{Z_{j}\widehat{\sigma}} )$. Hence
\begin{eqnarray*}
&& \bigl\|\mathbb{E}_{\mathbb{F}_m} \bigl[f (xS_{j}\eta _{j|m}
) - f (xZ_{j}\widehat{\sigma}\eta_{j|m} ) \bigr]
\bigr\|_1 \leq2 \bigl\||x|^p|\eta_{j|m}|^p
\bigr\|_1 \zeta_{p} (P_{S_{j}}, P_{Z_{j}\widehat{\sigma}} ).
\end{eqnarray*}
By Lemma~\ref{lemsigexpansion}, we have $ \|\eta_{j|m} \|
_p < \infty$. Note that $\mathbb{E} [S_{j} ] = \mathbb
{E} [Z_j ]
= 0$ and $ \|S_{j} \|_2^2 =  \|Z_j\widehat{\sigma
} \|_2^2$. Hence due to (\ref{eqlemapproxchar1}) and (\ref
{eqlemapproxchar2}), we may apply Theorem~3.1(1) in
\cite{dedeckerminimaldistance2009}, which gives us $\zeta_{p}
(P_{S_{j}}, P_{Z_{j}\widehat{\sigma}} ) \lesssim m^{-p/2 + 1}$. Hence
\begin{eqnarray*}
&&\bigl\|\mathbb{E}_{\mathbb{F}_m} \bigl[f (xS_{j}\eta _{j|m} )
- f (xZ_{j}\widehat{\sigma}\eta_{j|m} ) \bigr]
\bigr\|_1 \lesssim |x|^p m^{-p/2+1}.
\end{eqnarray*}
%
\item[\textit{Case} $p = 3$:] We may proceed as in the previous case, with the
exception that here we need to apply Theorem~3.2 in \cite
{dedeckerminimaldistance2009}. We may do so due to (\ref
{eqlemapproxchar1}) and (\ref{eqlemapproxchar2}). Hence we obtain
%
\begin{eqnarray}
&&\bigl\|\mathbb{E}_{\mathbb{F}_m} \bigl[f (xS_j\eta_{j|m} ) -
f (xZ_{j}\widehat{\sigma}\eta_{j|m} ) \bigr]
\bigr\|_1 \lesssim |x|^3 m^{-1/2}.
\end{eqnarray}
\end{longlist}
\upqed\end{pf}

Using Lemmas~\ref{lemtrunccharcosandsin} and~\ref
{lemapproxchar}, the triangle inequality gives the following
corollary, which proves (\ref{eqphiboundcor}).

\begin{cor}\label{corzolo}
Grant Assumption~\ref{assdependence}. Then
\begin{eqnarray*}
&&\bigl\|\varphi_i(t) -e^{-\sigma_{i|m}^2 t^2/2} \bigr\|_1 \lesssim
|t|^p m^{-p/2 + 1},\qquad i \in\{1,N\}.
\end{eqnarray*}
\end{cor}

\subsection{Some auxiliary lemmas}\label{secauxresults}

We will frequently use the following lemma, which is essentially a
restatement of Theorem~1 in \cite{sipwu}, adapted to our setting.

\begin{lem}\label{lemwuoriginal}
Put $p' = \min\{p,2\}$, $p \geq 1$. If $\sum_{l = 1}^{\infty}\sup_{k \in
\mathbb{Z}
} \|X_k-X_k^{(l,\prime)} \|_p < \infty$, then\vspace*{-3pt}
\begin{eqnarray*}
&& \|X_1 + \cdots+ X_n \|_{p} \lesssim
n^{1/p'}.
\end{eqnarray*}
For the sake of completeness, we sate this result in the general,
nontime-homogenous but stationary Bernoulli-shift context.
\end{lem}

Recall that\vspace*{-3pt}
\begin{eqnarray*}
 Y_j &=& \frac{1}{\sqrt{2m}} Y_j^{(1)},
\qquad Y_j^{(1)} = U_j + R_j,\\[-3pt]
S_{j|m}^{(**)} &=& \frac{1}{\sqrt{2m}}\sum
_{k = (2j-2)m +
1}^{(2j-1)m} X_k^{((2j - 2)m, **)}.
\end{eqnarray*}

\begin{lem}\label{lemboundR1}
Grant Assumption~\ref{assdependence}. Then:
\begin{longlist}[(ii)]
\item[(i)] $ \|S_{j|m}^{(**)} - Y_j \|_p \lesssim
m^{-1/2} (1 + \|R_j\|_p) \lesssim m^{-1/2}$\hspace*{1.5pt} for $j = 1,\ldots,N$,
\item[(ii)] $ \|Y_j \|_p < \infty$ for $j = 1,\ldots,N$.
\end{longlist}
\end{lem}

\begin{pf}
Without loss of generally, we assume that $j = 1$ since $m \thicksim m'$.
(i) We have the decomposition
\begin{eqnarray*}
&& \sqrt{2m} \bigl\|S_{j|m}^{(**)}- Y_j \bigr\|_p \leq\sum
_{k = 1}^m \bigl\|X_k -
X_k^{(**)} \bigr\|_p + \sum
_{k = 1}^m \bigl\|\mathbb{E}_{\mathbb
{F}_m}
[X_k ] \bigr\|_p + \|R_j \|_p.
\end{eqnarray*}
We will deal with all three terms separately. The triangle inequality gives
\begin{eqnarray*}
&& \sum_{k = 1}^m \bigl\|X_k -
X_k^{(**)} \bigr\|_p \leq\sum
_{k =
1}^{\infty} k \bigl\|X_k -
X_k' \bigr\|_p < \infty.
\end{eqnarray*}
Next, note that $\mathbb{E} [X_k^{(**)} |\mathbb{F}_m ]
= \mathbb{E}
[X_k ]= 0$ for $1 \leq k \leq m$. Hence it follows via the Jensen
and triangle inequalities that
\begin{eqnarray*}
 \Biggl\|\sum_{k = 1}^m\mathbb{E}
[X_k |\mathbb{F}_m ] \Biggr\|_p &=& \Biggl\|\sum
_{k = 1}^m\mathbb{E} \bigl[X_k -
X_k^{(**)} |\mathbb {F}_m \bigr] \Biggr\|_p
\\
&\leq & \sum_{k = 1}^m \bigl\|X_k -
X_k^{(**)} \bigr\|_p \leq\sum
_{k = 1}^{\infty}k \bigl\|X_k - X_k'
\bigr\|_p < \infty.
\end{eqnarray*}
Similarly, since $X_k - \mathbb{E}_{\mathbb{F}_m} [X_k ]
\stackrel{d}{=}
\mathbb{E}_{\mathbb{F}_m} [X_k^{(k-m,*)} - X_k ]$ for $m+1
\le k \leq2m$,
we have
\begin{eqnarray*}
 \Biggl\|\sum_{k = m+1}^{2m} X_k -
\mathbb{E}_{\mathbb{F}_m} [X_k ] \Biggr\|_p &=& \Biggl\|\sum
_{k = m+1}^{2m} \mathbb{E}_{\mathbb
{F}_m}
\bigl[X_k^{(k-m,*)} - X_k \bigr] \Biggr\|_p
\\
&\leq & \sum_{k =
m+1}^{2m} \bigl\|
X_k^{(k-m,*)} - X_k \bigr\|_p \leq\sum
_{k =
1}^{\infty}k \bigl\|X_k'
- X_k \bigr\|_p < \infty.
\end{eqnarray*}
Combining all three bounds gives (i). This implies that for (ii), it
suffices to show that $ \|U_1 \|_p \lesssim\sqrt{m}$. Using
the above bounds and Lemma~\ref{lemwuoriginal}, we get
\begin{eqnarray*}
&&\|U_1 \|_p \leq \Biggl\|\sum_{k = 1}^m
X_k \Biggr\|_p +\sum_{k = 1}^m
\bigl\|\mathbb{E}_{\mathbb{F}_m}[X_k] \bigr\|_p \lesssim\sqrt{m}.
\end{eqnarray*}
\upqed\end{pf}

\begin{lem}\label{lemauxHsigcomp}
Grant Assumption~\ref{assdependence}, and let  $H_j$ and
$\mathcal{H}_j$, $j \in \mathcal{J}$ be as in (\ref
{eqthmauxcharequationbound2.91}). Then
\begin{eqnarray*}
&& P \bigl(\mathbb{E}_{\mathcal{H}_j} \bigl[H_j^2 \bigr]
\leq\widehat {\sigma }_{m-l}^2 \bigr) < 1/7.
\end{eqnarray*}
\end{lem}

\begin{pf}
Since $m \thicksim m'$, it suffices to treat the case $j = 1$. Recall
that $\mathcal{H}_1 = \mathcal{G}_1^{(l)} = \sigma(\mathcal{E}_l \cup \{\epsilon_{m+1},\ldots,\epsilon_{2m}\})$ and
\begin{eqnarray*}
\sqrt{m-l}H_1 &= & \sum_{k= l+1}^m
X_k - \mathbb{E}_{\mathcal{G}_1^{(l)}} [X_k ] - R_1
+ \mathbb{E}_{\mathcal{G}_1^{(l)}} [R_1 ],
\\
2(m-l) \widehat{\sigma}_{m-l}^2 &=& \mathbb{E}_{\mathcal{E}_l}
\Biggl[ \Biggl(\sum_{k = l+1}^m
X_k^{(k-l,*)} \Biggr)^2 \Biggr].
\end{eqnarray*}
Let $I_k^{(l)} = \mathbb{E}_{\mathcal{E}_l}[X_k] + (R_1 - \mathbb{E}_{\mathcal{G}_1^{(l)}}[R_1])\mathbh{1}(k = l+1)$. Using
$a^2-b^2 = (a-b)(a+b)$ and applying the Cauchy--Schwarz and Jensen
inequalities then yields
\begin{eqnarray*}
&&(m-l) \bigl\|\mathbb{E}_{\mathcal{H}_1} \bigl[H_1^2 \bigr]
- 2\widehat {\sigma}_{m-l}^2 \bigr\|_1
\\
&&\qquad = \Biggl\|\mathbb{E}_{\mathcal{E}_l} \Biggl[ \Biggl(\sum_{k =
l+1}^m
X_k - X_k^{(k-l,*)} - I_k^{(l)} \Biggr)
\Biggl(\sum_{k = l+1}^m X_k +
X_k^{(k-l,*)} - I_k^{(l)} \Biggr) \Biggr] \Biggr\|_1
\\
&&\qquad\leq \Biggl\|\sum_{k = l+1}^m X_k -
X_k^{(k-l,*)} - I_k^{(l)} \Biggr\|_2 \Biggl\|\sum
_{k =
l+1}^m X_k + X_k^{(k-l,*)}
- I_k^{(l)} \Biggr\|_2
\\
&&\qquad\stackrel{\mathrm{def}} {=}  I_1(l,m) \mathit{II}_2(l,m).
\end{eqnarray*}
By\vspace*{1pt} Lemma~\ref{lemboundR1} and the arguments therein, it follows that
$I_1(l,m) = \mathcal{O} (1 )$, uniformly for $0<l<m$.
Similarly, one
obtains that $(m-l)^{-1/2}\mathit{II}_2(l,m) = \mathcal{O} (1 )$. Hence
%
\begin{eqnarray}
&& \bigl\|\mathbb{E}_{\mathcal{H}_j} \bigl[H_j^2 \bigr] - 2
\widehat {\sigma }_{m-l}^2 \bigr\|_1 \lesssim
\frac{1}{\sqrt{m-l}}.
\end{eqnarray}
We then have that
\begin{eqnarray*}
P \bigl(\mathbb{E}_{\mathcal{H}_j} \bigl[H_j^2 \bigr]
\leq\widehat {\sigma }_{m-l}^2 \bigr) &\leq &  P \bigl( \bigl|
\mathbb{E}_{\mathcal
{H}_j} \bigl[H_j^2 \bigr] - 2
\widehat{\sigma}_{m-l}^2 \bigr| \geq \widehat{\sigma
}_{m-l}^2 \bigr)
\\
&\leq & \widehat{\sigma}_{m-l}^{-2} \bigl\|\mathbb{E}
_{\mathcal{H}_j} \bigl[H_j^2 \bigr] - 2\widehat{
\sigma}_{m-l}^2 \bigr\|_1 \lesssim\frac{1}{\sqrt{m-l}}
\end{eqnarray*}
by Markov's inequality. Hence the claim follows if $m - l$ is large
enough. Note that more detailed computations, as in Lemma~\ref
{lemsigexpansion}, would give a more precise result. However, the
current version is sufficient for our needs.
\end{pf}

%


\section*{Acknowledgments}

I would like to thank the Associate Editor and the anonymous reviewer
for a careful reading of the manuscript and the comments and remarks
that helped to improve and clarify the presentation. I also thank Istvan Berkes and Wei Biao Wu for stimulating discussions.
Special thanks to Florence Merlev\`{e}de for pointing out a few errors.

%





\printaddresses
\end{document}